\documentclass[a4paper, final, 10pt]{article}

\usepackage{amsmath,amsfonts,amssymb,amstext}
\usepackage{bm}
\usepackage{graphicx}
\usepackage{hyperref}
\usepackage{subfig}

\usepackage[utf8]{inputenc}
\usepackage[english]{babel}

\usepackage{pgfplots}
\usepackage{verbatim}
\usepackage{bm}

\usepackage{tikz}
\usetikzlibrary{matrix,calc}
\usepackage{algorithm,algorithmic}

\newcommand{\grad}{\mathop{\rm \nabla }\nolimits}

\renewcommand{\d}[2]{\frac{d #1}{d #2}} 
\newcommand{\pd}[2]{\frac{\partial #1}{\partial #2}}

\newcommand{\pdn}[1]{\frac{\partial #1}{\partial n}} 

\graphicspath{ {./figs/}}

\begin{document}

%

\title{Numerical recovery of the piecewise constant leading coefficient of an elliptic equation\thanks{ The work was supported by the grant of the President of the Russian Federation for state support of young scientists MK-1131.2020.1, Mega grant of Russian Government N14.Y26.31.0013, and RFBR grant N20-01-00207. }}
\author{Aleksandr E. Kolesov \thanks{The corresponding author} \thanks{North-Eastern Federal University, 48 Kulakovskogo str., Yakutsk 677000, Russia},  
Petr N. Vabishchevich \footnotemark[2] \thanks{Nuclear Safety Institute of RAS, 52 B. Tulskaya str., Moscow 115191, Russia}
}
\date{}

\maketitle

\begin{abstract}
We propose a numerical algorithm for the reconstruction of a piecewise constant leading coefficient of an elliptic problem. 
The inverse problem is reduced to a shape reconstruction problem. 
The proposed algorithm is based on the minimization of a cost functional where a control function is the right-hand side of an auxiliary elliptic equation for a level set representation of unknown shape. 
The numerical implementation is based on the finite element method and the open-source computing platform \texttt{FEniCS}.
The performance of the algorithm is demonstrated on computationally simulated data.
\end{abstract}

\textbf{Keywords}: coefficient inverse problem, elliptic equation,  level set,  adjoint method,  finite element method

%
%
%
%


\section{Introduction}

The coefficient inverse problem of identifying the unknown leading coefficient of an elliptic equation is the basis of electrical impedance tomography (EIT) and electrical resistivity tomography (ERT) non-invasive technique for investigation of the internal structure of bodies from voltage and current boundary measurements. The inverse problem for elliptic equation  has numerous applications in geophysics \cite{Chambers1999, Chambers2006}, medical imaging \cite{Bodenstein2009a, Putensen2019} and nondestructive testing \cite{McCarter1989, Karhunen2010}.

The coefficient inverse problem for elliptic equation  was first posed in \cite{Calderon1980}. It is known to be nonlinear and ill-posed and, therefore, special algorithms are needed to solve such problems numerically. Also, the issues of stability and uniqueness are crucial for the development of robust numerical algorithms. The uniqueness of a large class of isotropic coefficients was considered in many works \cite{Sylvester1986, Novikov1988, Nachman1996, Brown1997, Astala2006}. In \cite{Astala2006}, the optimal regularity condition for coefficient in two dimensions was obtained. The stability of the inverse problem was investigated in \cite{Alessandrini1988, Allers1991, Barcelo2001a}. 
To achieve stability, it is required to place some constraints on the coefficient. 
For example, in \cite{Alessandrini1988}, a logarithmic stability estimate was obtained.

There are many reconstruction methods and procedures for solving the inverse problem \cite{Cheney1999, Borcea2002}. These methods can be divided into two main groups: non-iterative and iterative methods. 
The non-iterative methods include  factorization \cite{Bruhl2000, Gebauer2007},  layer stripping \cite{Somersalo1991, Sylvester1992},  D-bar methods \cite{Siltanen2001, Knudsen2009}, NOSER \cite{Cheney1990}.
The factorization method is a method based on an explicit criterion for detecting inhomogeneities inside bodies \cite{Bruhl2001}. 
In layer-stripping algorithms, one needs first to find the unknown coefficient on the boundary of the body and then progress inside layer by layer. 
The D-bar method is based on evaluating a nonlinear Fourier transform of the coefficient from EIT data and the inversion of the transform. 
NOSER (Newton's One-Step Error Reconstructor) is an algorithm based on the minimization of an error functional, but it takes only one step of Newton's method with constant coefficient as an initial guess. Most of the calculations, including the gradient of the functional, can be done analytically.

Iterative methods for solving the inverse problems include minimization algorithms based on either least squares \cite{Yorkey1987, Dobson1992, Borcea2001} or equation-error \cite{Wexler1985, Kohn1990} formulations. 
Least squares methods usually iteratively minimize the norm of the difference between electrical potential due to the applied current and the measured potential on the boundary. 
Equation-error approaches are also known as variational methods, derived from Dirichlet and Thompson variational principles \cite{Borcea2002}. 
The convergence of iterative methods can be ensured using some regularization technique, which typically means that a specific regularization term is added to a functional. 
Iterative methods require calculating first and, in some optimization methods, second derivatives of the objective functionals. 
The adjoint method can be used to calculate both derivatives \cite{Dorn2000} efficiently.  

In many applications of the coefficient inverse problem, one can assume that objects under investigation contain several materials with piecewise constant coefficients 
This assumption allows us to reduce the inverse problem of recovering the distribution of  coefficient inside the body to the problem of reconstructing shapes of the materials. 
Among the shape reconstruction methods, the level set method is known to be the most powerful one. 

The level set method was first proposed in \cite{Osher1988} to track evolving interfaces.
There, an evolving domain is represented by a  continuous level set function. Then, the motion of this domain is expressed via a Hamilton-Jacobi equation for the level set function. 
In \cite{Santosa1996}, this approach was first used to solve inverse problems where the desired unknown is a  characteristic function of some geometry. 
Next, the unknown is represented by a level set function.
The evolution of the level set function minimizes a functional and leads to a solution of the inverse problem.  

 In \cite{Osher2001}, the level set method was used to solve the inverse problem associated with shape optimization for the eigenvalue problem for the Laplace equation. This approach was used to solve the inverse problem for electrical impedance tomography in \cite{Ito2001, Chan2004, Chung2005}. 
More recently, the level set based methods were applied to elliptic inverse problems in \cite{Agnelli2018, Lin2018, Liu2019}. Surveys on level set methods for solving inverse problems can be found in \cite{Tai2004, Burger2005a, Gibou2018}.

In this paper, we present a numerical algorithm based on the level set idea for solving the inverse problem associated with electrical impedance tomography. 
We assume that the leading coefficient is a piecewise constant function, and the values of coefficients are known. 
The main idea is to implicitly represent the interface between regions with known coefficients as the zero value of an auxiliary elliptic equation's solution.
The algorithm is based on minimizing the squared norm of the difference between potentials due to applied currents and measured potentials. We previously used this approach to successfully recover a piecewise constant lower coefficient  \cite{Kolesov2019} and the right-hand side \cite{Ivanov2019} of an elliptic equation.

The paper is organized as follows. In the next section, we introduce the inverse problem. The details of the proposed reconstruction algorithm, the cost functional, and the calculation of its derivative using the adjoint method are discussed in section 3. Then, we present a series of numerical experiments to show our algorithm's ability to recover the unknown interface from noisy data. 
The final section is the conclusion.

\section{Problem statement}

Let $\Omega \subset \mathbb{R}^2=d$ ($d=2,3$) be the bounded domain  with sufficiently smooth boundary $\partial\Omega$ and $\sigma(\bm x) \in L^{\infty}(\Omega)$ is the coefficient such that $\sigma(\bm x) \geq \sigma_0 > 0$.
We consider the boundary value problem for elliptic equation
\begin{equation} 
\label{eq:forward}
- \nabla \cdot \sigma (\bm x) \nabla u = 0, \quad \bm x \in \Omega,
\end{equation}
\begin{equation}
\label{eq:g}
\sigma (\bm x)  \pdn{u} = g (\bm x), \quad \bm x \in \partial\Omega,
\end{equation}
where $u(\bm x)$ is the electrical potential, $\bm n$ is  the unit outward normal to $\partial \Omega$ and $g(\bm x)$ is the applied current density. In addition, both $u(\bm x)$ and  $g(\bm x)$ must satisfy the following constraints  
\begin{equation}
\label{eq:uds}
\int_{\partial \Omega} u(\bm x) ds = 0,
\end{equation}
\begin{equation}
\label{eq:gds}
\int_{\partial \Omega} g(\bm x) ds = 0.
\end{equation}
The problem (\ref{eq:forward})--(\ref{eq:gds}) correspond to the continuum model for electrical impedance tomography. 
The inverse problem is to determine the distribution of  coefficient $\sigma(\bm x)$ inside $\Omega$ using a set of given values of applied current density $g(\bm x)$ on $\partial \Omega$ and the corresponding measured values of potential $u( \bm x)$ on $\partial \Omega$ . The set of $g(\bm x)$ and $u(\bm x)$ on $\partial \Omega$  is also known as the Neumann-to-Dirichlet or current-to-voltage map in problems of electrical impedance tomography.  It is well known that this inverse problem does not have unique solution. Therefore, we need to narrow the class of admissible solutions.  
 


In many applications, such as medical imaging, geophysics, and nondestructive testing, it is a-priori known that the object to be imaged contains several materials with piecewise constant conductivities. Let $N$ be the number of materials, $D_i$ be the subdomain containing material with  $\sigma_i$, ($i = 1, \dots, N$),  $\Omega = \cup_{i=1}^N D_i$.  Then,  the distribution of  coefficient $\sigma (\bm x)$ can be represented as
\begin{equation*} 
\sigma (\bm x) = \sum_{i=1}^N  \sigma_i \chi_{D_i}(\bm x),
\end{equation*}
where $\chi_{D_i}(\bm x)$  ($i=1,\dots, N$) is the characteristic function of the subdomain $D_i$:
\begin{equation}\label{eq:chi}
\chi_i(\bm x) = \begin{cases}
1, & \bm x \in D_i,\\
0, & \bm x \in \Omega \backslash D_i.
\end{cases}
\end{equation}
In this work, for simplicity  we restrict ourselves with only two materials with conductivities $1$ and $2$ an denote by $D$ the subdomain with coefficient $2$. Therefore, the coefficient distribution $\sigma (\bm x)$ can be written as
\begin{equation} \label{eq:coef}
\sigma (\bm x) = 1 + \chi_{D}(\bm x)),
\end{equation}
with $\chi_{D} (\bm x)$ is the characteristic function of the subdomain $D$. 

In this case, the inverse problem of determining the distribution of coefficient $\sigma (\bm x)$ reduces to the reconstruction of the shape of the subdomain $D$.

\section{Reconstruction algorithm}

In this section, we propose a new algorithm for reconstruction of the coefficient $\sigma (\bm x)$.
The algorithm is based on minimizing  a cost functional  using a gradient method.

\subsection{Cost functional}

We introduce a level set function $q(\bm x)$, which describes subdomain  $D$ as follows
\begin{equation}\label{eq:levelset}
\begin{cases}
q(\bm x) \geq 0, & \bm x \in D,\\
q(\bm x) < 0, & \bm x \in \Omega \backslash D,
\end{cases}
\end{equation}
and the Heaviside function $H(q)$ 
\begin{equation}
\label{eq:heaviside}
H(q(\bm x)) = 
\begin{cases}
1, & q(\bm x) \geq 0, \\
0, & q(\bm x) < 0.
\end{cases}
\end{equation}
Then, the coefficient $\sigma (\bm x)$ (\ref{eq:coef}) can be defined as
\begin{equation}\label{eq:c_new}
\sigma (\bm x) = 1 + H(q(\bm x)).
\end{equation}
Clearly, to determine  $\sigma (\bm x)$  it is sufficient to identify the level set function $q(\bm x)$. Many different level set functions were used to solve the inverse problem of electrical impedance tomography \cite{Dorn2006}. Commonly, the level set function can be determined using a signed distance function, which can be updated by solving the Hamilton-Jacobi equation.

In this work, the level set function is the solution of the following elliptic equation:
\begin{equation} \label{eq:aux}
- \gamma \Delta q + q = f(\bm x), \quad \bm x \in \Omega,
\end{equation}
\begin{equation} \label{eq:aux_bc}
q = 0, \quad \bm x \in \partial\Omega,
\end{equation}
where $\gamma = \mbox{const} > 0$ is the parameter.  The key of our approach is to determine the right hand side $f(\bm x)$ such as the solution $q(\bm x)$ of problem (\ref{eq:aux}), (\ref{eq:aux_bc}) describes the desired coefficient $\sigma (\bm x)$ (\ref{eq:c_new}).
In fact, the parameter $\gamma$ can be seen as a smoothing parameter for function $f(\bm x)$ and $q(\bm x)$ is smoothed out function. For $\gamma = 0$, we have $q(\bm x) = f(\bm x)$.  Note that we can employ the different differential operator instead of $-\Delta$.

Let $M$ be the number of measurements. For $1 \leq j \leq M$, let $g_j(\bm x)$ be a given applied current density on $\partial \Omega$ and $m_j(\bm x)$ be the corresponding measurement of the potential on $\partial \Omega$.  
To find $f(\bm x)$, we minimize the following least-squares cost functional
\begin{equation} \label{eq:J}
J(f) = \frac{1}{2} \sum_{j=1}^M \int_{\partial \Omega} \left|u_j(\bm x; f) - m_j(\bm x) \right|^2 \, \mathrm{d}s.
\end{equation} 
where $u_j(\bm x; f)$, $j = 1, \dots, M$ are the solutions of the problems
\[
- \nabla \cdot \sigma (\bm x) \nabla u_j = 0, \quad \bm x \in \Omega,
\]
\[
\sigma(\bm x)  \pdn{u_j} = g_j (\bm x), \quad \bm x \in \partial\Omega.
\]
Here, $\sigma(\bm x) $ is a coefficient corresponding to $f(\bm x)$ via (\ref{eq:c_new})--(\ref{eq:aux_bc}).  
Note that in (\ref{eq:J}) the functions $m_j(\bm x)$ correspond to the solution of above problem for the desired coefficient.
Both $u_j(\bm x; f)$ and $m_j(\bm x)$ must satisfy (\ref{eq:uds}).

\subsection{Variational formulations}

For discretization in space we use the finite element method, so we need to obtain variational forms of boundary value problems  (\ref{eq:forward})--(\ref{eq:gds}) and   (\ref{eq:aux}), (\ref{eq:aux_bc}).   
First, we define the functional spaces:
\[
V = \left\{v \in H^1(\Omega): \int_{\partial \Omega} v \, d\bm s = 0\right \},  \quad Q = \left\{v \in H^1(\Omega): v(\bm x ) = 0,  \, \bm x \in \partial \Omega  \right \}, 
\]
where $H^1(\Omega)$ is Sobolev space. 
We multiply equation (\ref{eq:forward}) by a test function $v \in V$, integrate the resulting equation over $\Omega$ and perform integration by parts to eliminate second-order derivative of $u$:
\begin{equation*}
\int_{\Omega} \sigma \nabla u \,  \nabla v \, d\bm x - \int_{\partial \Omega} \sigma \pdn{u} \, v \, d\bm s = 0.
\end{equation*}
Taking into account the boundary condition (\ref{eq:g}), yields
\begin{equation*}
\int_{\Omega} \sigma \nabla u \,  \nabla v \, d\bm x = \int_{\partial \Omega} g \, v \, d\bm s.
\end{equation*}
Next, we define the following bilinear form
\begin{equation*}
a(u,v) = \int_{\Omega} \sigma \nabla u \,  \nabla v \, d\bm x, 
\end{equation*}
and linear form
\begin{equation*}
L_a(v) = \int_{\partial \Omega} g \, v \, d\bm s.
\end{equation*}
Then the variational formulation of problem (\ref{eq:forward})--(\ref{eq:gds}) read as: find $u \in V$ such as
\begin{equation}
\label{eq:forward_var}
a(u,v) = L_a(v),  \quad \forall v \in V.
\end{equation}

The variational formulation of problem  (\ref{eq:aux}), (\ref{eq:aux_bc}) is derived similarly: find $q \in Q$ such that  
\begin{equation}
\label{eq:aux_var}
b(q, w) = L_b(w),  \quad \forall w \in Q,
\end{equation}
where
\begin{equation*}
b(q, w) = \int_{\Omega} \gamma \nabla q \,  \nabla w \, d\bm x +  \int_{\Omega}  q \, w \, d\bm x,
\end{equation*}
\begin{equation*}
L_b(w) = \int_{\Omega}  f \, w \, d\bm x,
\end{equation*}

Note that both bilinear forms $a(u,v) $ and $b(q, w)$ are symmetric
\begin{equation*}
a(u,v) = a(v,u), \quad b(q, w) = b(w, q).
\end{equation*}
This property is useful for calculating the gradient of $J(f)$.

\subsection{Gradient of functional}

To minimize the functional $J(f)$ we use a gradient based method.Thus, we need to calculate the gradient of $J(f)$ with respect to the function $f(\bm x)$. First, note that  $J(f)$ is the functional of  the functions $u_j(\bm x)$, $j = 1, \dots, M,$ which depend on  coefficient $\sigma(\bm x)$. In turn, $\sigma(\bm x)$ is the function of the level set $q(\bm x)$ via (\ref{eq:c_new}). Finally, the function $q(\bm x)$  also depends on the objective function $f(\bm x)$ by (\ref{eq:aux}). Therefore, by the chain rule, to compute the gradient of $J$ we use the following equation 
\begin{equation}
\label{eq:dJdf1}
\d{J}{f} = \sum_{j=1}^{M} \pd{J}{u_j} \pd{u_j}{\sigma} \pd{\sigma}{q} \pd{q}{f}.
\end{equation}
The terms $\partial J / \partial u_j$ and $\partial \sigma / \partial q$ are straightforward to compute using the following equations:
\begin{equation}
\label{eq:djdu}
\pd{J}{u_j}  = u_j - m_j,
\end{equation}
and
\begin{equation}
\label{eq:dsigmadq}
\pd{\sigma}{q} = \delta(q(\bm x)),
\end{equation}
where $\delta (\cdot)$ is the Dirac delta function.
 By contrast, $\partial u_j / \partial \sigma$ and $\partial q / \partial f$ are rather difficult to compute. 
So, we use the adjoint method to compute the gradient of $dJ/df$.
 
First, by taking the derivative of (\ref{eq:forward_var}) with respect to $q$, we get
\[
\pd{a}{u_j} \pd{u_j}{\sigma}  + \pd{a}{\sigma} = 0, \quad j = 1, \dots, M.
\] 
Since $\partial a / \partial u_j$ is invertible, the following equation for $\partial u_j/\partial\sigma$ is obtained
\begin{equation}
\label{eq:dudsigma}
\pd{u_j}{\sigma} = - \left( \pd{a}{u_j} \right)^{-1} \pd{a}{\sigma} .
\end{equation}
By analog, we can derive the equation for $\partial q/\partial f$ from (\ref{eq:aux_var})
\begin{equation}
\label{eq:dqdf}
 \pd{q}{f}  = \left( \pd{b}{q} \right)^{-1} \pd{L_b}{f} 
\end{equation}
Substituting  (\ref{eq:dudsigma}) and (\ref{eq:dqdf}) into (\ref{eq:dJdf1}) yields 
\begin{equation*}
\d{J}{f} =- \sum_{j=1}^{M} \pd{J}{u_j} \left( \pd{a}{u_j} \right)^{-1} \pd{a}{\sigma}  \pd{\sigma}{q} \left( \pd{b}{q} \right)^{-1} \pd{L_b}{f}.
\end{equation*}

Now we take adjoint of the above equation
\begin{equation}
\label{eq:dJdf2}
\d{J}{f}^*=\pd{L_b}{f}^* \left( \pd{b}{q} \right)^{-1*}   \pd{\sigma}{q}^* \sum_{j=1}^{M} \pd{a}{\sigma}^*  \left(  \pd{a}{u_j} \right)^{-1*} \d{J}{u_j}^* 
\end{equation}
Next, we define new variables $z_j$ as follows
\[
z_j = \left( \pd{a}{u_j} \right)^{-1*} \, \d{J}{u_j}^*, \quad j = 1, \dots, M,
\]
and obtain the following adjoint equations associated with (\ref{eq:forward_var})
\begin{equation*}
\pd{a}{u_j}^{*} \, z_j = \d{J}{u_j}^*, \quad j = 1, \dots, M.
\end{equation*}
Taking into account that the bilinear form $a$ is symmetric and  (\ref{eq:djdu}) we get 
\begin{equation}
\label{eq:forward_adjoint}
a(z_j,v) = \int_{\partial \Omega} (u_j - m_j) \, v \, ds, \quad j = 1, \dots, M.
\end{equation}
Substituting the adjoint solutions $z_j$ into (\ref{eq:dJdf2}) leads to
\begin{equation}
\label{eq:dJdf3}
\d{J}{f}^*=\pd{L_b}{f}^* \left( \pd{b}{q} \right)^{-1*}   \pd{\sigma}{q}^* \sum_{j=1}^{M} \pd{a}{\sigma}^*   z_j 
\end{equation}

Next, we define another adjoint variable $\lambda$ as follows
\begin{equation*}
\lambda = \left( \pd{b}{q} \right)^{-1*}   \pd{\sigma}{q}^* \sum_{j=1}^{M} \pd{a}{\sigma}^*   z_j 
\end{equation*}
and the adjoint equation associated with (\ref{eq:aux_var})
\begin{equation*}
\pd{b}{q}^{*}  \lambda = \pd{\sigma}{q}^* \sum_{j=1}^{M} \pd{a}{\sigma}^*   z_j 
\end{equation*}
The bilinear form $b$ is also symmetric, so we have the adjoint problem associated with (\ref{eq:aux_var})
\begin{equation}
\label{eq:aux_adjoint}
b(\lambda, w) =  \delta(q) \sum_{j=1}^M \int_{\Omega} \grad u_j \cdot \grad z_j \, w \, dx
\end{equation}

Finally, we substitute the adjoint solution $\lambda$ into (\ref{eq:dJdf3})
\begin{equation*}
\d{J}{f}^* = \pd{L_b}{f}^* \, \lambda
\end{equation*}
which yields the final equation for computing the gradient of functional $J$ with respect to $f$: 
\begin{equation}
\label{eq:dJdf4}
\d{J}{f} = \int_{\Omega} \lambda \, w \, dx.
\end{equation}
To compute the gradient $dJ/df$ (\ref{eq:dJdf4})  we need to solve the adjoint problems (\ref{eq:forward_adjoint}) and (\ref{eq:aux_adjoint}). 

\subsection{Algorithm}

In this section, we propose the new algorithm for solving our coefficient inverse problem. It is based on finding $f(\bm x)$ such that the solution  $q(\bm x)$ of problem (\ref{eq:aux}), (\ref{eq:aux_bc}) describe the desired coefficient  $\sigma (\bm x)$ via (\ref{eq:c_new}). To find $f(\bm x)$ we minimize the functional $J(f)$(\ref{eq:J}) using a gradient based iterative procedure:

\begin{algorithm}
\caption*{\textbf{Reconstruction algorithm}}
\begin{algorithmic}[1]
\STATE{
\textbf{input} $M$, $g_j$, $m_j$, $j=1,\dots,M$, $f^0$, $ \gamma, \alpha, tol, K$ \COMMENT{$K$ is a large number}
}
\STATE{
\textbf{output} $\sigma(\bm x)$
}
\FOR{$k=0$ \TO $K$}
\STATE{
$q^k(\bm x) \leftarrow$ the solution of  (\ref{eq:aux_var})
}
\STATE{
$\sigma^k(\bm x) \leftarrow 1 + H(q^k(\bm x))$
}
\FOR{$j=1$ \TO M}
\STATE{
$u_j^k \leftarrow$ the solution of  (\ref{eq:forward_var}) for given $g_j(\bm x)$
}
\ENDFOR
\STATE{
$J^k(f^k) \leftarrow \frac{1}{2} {\displaystyle\sum_{j=1}^M \int_{\partial \Omega} \left|u_j(f;\bm x) - m_j(\bm x) \right|^2 \, \mathrm{d}s }$
}
\FOR{$j=1$ \TO M}
\STATE{
$z_j^k(\bm x) \leftarrow$ the solution of adjoint problem  (\ref{eq:forward_adjoint}) for given $m_j(\bm x)$
}
\ENDFOR
\STATE{
$\lambda^k(\bm x) \leftarrow$ the solution of adjoint problem  (\ref{eq:aux_adjoint})
}
\STATE{
$\cfrac{dJ^k}{df^k}  \leftarrow  \int_{\Omega} \lambda \, w \, dx.$
}

\IF{ $\left\|\cfrac{dJ^k}{df^k} \right \|_{L^{\infty}(\Omega)}  < tol$}
\STATE{
$\sigma(\bm x) \leftarrow \sigma^k(\bm x)$
}
\STATE{\textbf{stop}}
\ENDIF
\STATE{
$f^{k+1} \leftarrow f^k - \beta^k \cfrac{dJ^k}{df^k} $
}
\COMMENT{step size $\beta^k $ is calculated by line search}
\ENDFOR
\STATE{
$\sigma(\bm x) \leftarrow \sigma^k(\bm x)$
}
\end{algorithmic}
\end{algorithm}

The stopping criteria $tol$ will be discussed later. The numerical implementation of the reconstruction is performed in Python. For solving partial differential equations, we use an open-source computing platform FEniCS \cite{AlnaesBlechta2015a}. 

\section{Numerical experiments}

In this section, we examine the effectiveness of the proposed algorithm. A series of numerical experiments have been performed to confirm that our algorithm can recover the interface between two materials with piecewise constant coefficients from computationally simulated data. 

\subsection{Implementation details}

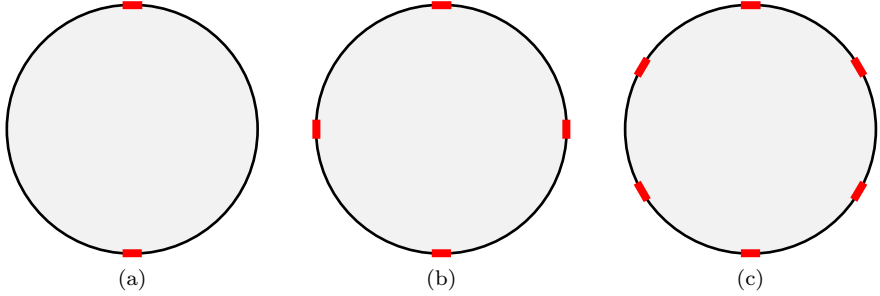
\begin{figure}
\begin{center}
\subfloat[\label{fig:2elecrodes}]{
\begin{tikzpicture}[scale=0.55]
\draw [line width=1.0,fill=black!5!white] (0,0) circle (3);
\draw [line width=3.0, red] (-0.23,2.99) -- (0.24,2.99);
\draw [line width=3.0, red] (0.23,-2.99) -- (-0.23,-2.99);
\end{tikzpicture}
}\hspace{1em}
\subfloat[\label{fig:4elecrodes}]{
\begin{tikzpicture}[scale=0.55]
\draw [line width=1.0,fill=black!5!white] (0,0) circle (3);
\draw [line width=3.0, red] (-0.23,2.99) -- (0.24,2.99);
\draw [line width=3.0, red] (2.99,0.23) -- (2.99,-0.23);
\draw [line width=3.0, red] (0.23,-2.99) -- (-0.23,-2.99);
\draw [line width=3.0, red] (-2.99,0.23) -- (-2.99,-0.23);
\end{tikzpicture}
}\hspace{1em}
\subfloat[\label{fig:6elecrodes}]{
\begin{tikzpicture}[scale=0.55]
\draw [line width=1.0,fill=black!5!white] (0,0) circle (3);
\draw [line width=3.0, red] (-0.23,2.99) -- (0.24,2.99);
\draw [line width=3.0, red] (2.47,1.7) -- (2.71,1.29);
\draw [line width=3.0, red] (2.71,-1.29) -- (2.47,-1.7);
\draw [line width=3.0, red] (0.23,-2.99) -- (-0.23,-2.99);
\draw [line width=3.0, red] (-2.47,-1.71) -- (-2.71,-1.29);
\draw [line width=3.0, red] (-2.71,1.29) -- (-2.47,1.7);
\end{tikzpicture}
}
%
\end{center}
\caption{Distribution of equidistant electrodes (in red) along the boundary of the domain: a) 2 electrodes, b) 4 electrodes, c) 6 electrodes}
\label{fig:electrodes}
\end{figure}

For numerical computation, we consider unit circle domain $\Omega$. 
To generate boundary data for potential $m$ we employ sets of equidistant electrodes along the boundary with width $\theta = \pi / 20$. See Figure \ref{fig:electrodes} for a scheme of distribution of 2, 4, and 6 electrodes. The following current injection pattern is used for measurements: we set $g(\bm x) = 1$ on one electrode, $g(\bm x) = -1$ on the opposite electrode, and $g(\bm x) = 0$ elsewhere.
 In this manner, we obtain one independent measurement for two electrodes, two measurements for four electrodes, etc. 

\begin{figure}
\begin{center}
\subfloat[\label{fig:mesh_ellipse}]{\includegraphics[width=0.33\textwidth]{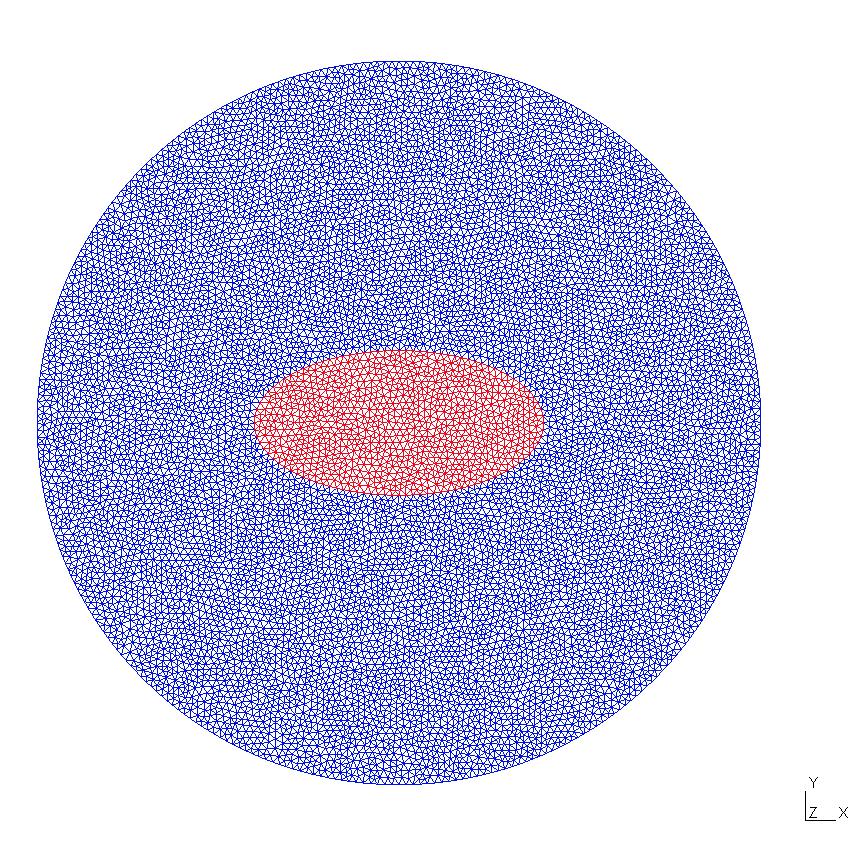}} 
\subfloat[\label{fig:mesh_circles}]{\includegraphics[width=0.33\textwidth]{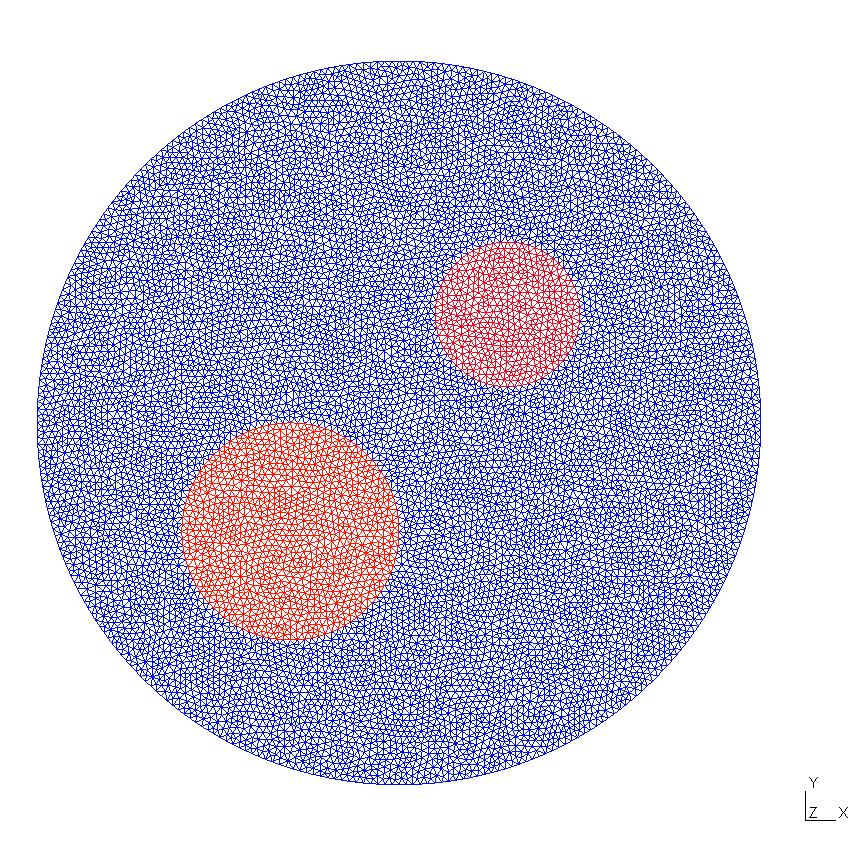}} 
\subfloat[\label{fig:mesh}]{\includegraphics[width=0.33\textwidth]{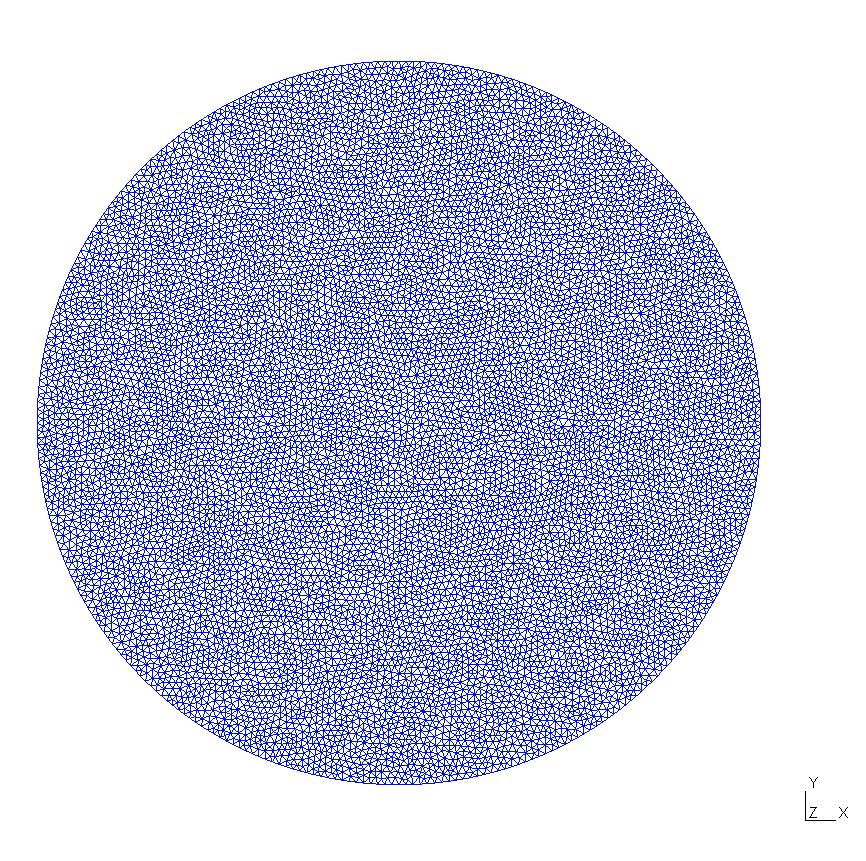}} 
\end{center}
\caption{Computational meshes used for: (a) data generation for ellipse,  (b) data generation for circles,  (d) reconstruction}
\label{fig:meshes}
\end{figure}

In our experiments, we reconstruct two images: ellipse and two circles with different diameters. 
To generate boundary data $m_j(\bm x)$, $j = 1, \dots, M$, $\bm x \in \partial \Omega$ for each test case we solve the forward problem (\ref{eq:forward})--(\ref{eq:gds}) using an unstructured mesh adapted to represent internal structure of imaged materials. For reconstruction purposes, we use a different mesh with no knowledge of imaged materials. All meshes are created using the open source software for mesh generation \texttt{Gmsh}. 
The computational meshes used for data generation and reconstruction are shown in Figure \ref{fig:meshes}. Note that the mesh in Figure \ref{fig:mesh} is used to reconstruct the images. 

For numerical solution we need to approximate the Heaviside function $H(q(\bm x))$ (\ref{eq:heaviside}) and the delta function $\delta (q(\bm x))$ by the following smooth functions $H_{\alpha} (q(\bm x))$ and $\delta_{\alpha} (q(\bm x))$
\begin{equation}
\label{eq:heaviside_alpha}
H_{\alpha}(q(\bm x)) = 
\begin{cases}
0, & q(\bm x) < 0, \\
\displaystyle{\frac{1}{2} - \frac{1}{2} \cos \left( \pi \frac{ q}{\alpha} \right)}, & 0 \leq q(\bm x) < \alpha \\
1, & q(\bm x) \geq \alpha,
\end{cases}
\end{equation} 
and
\begin{equation}
\label{eq:delta_alpha}
\delta_{\alpha}(q(\bm x)) = 
\begin{cases}
0, & q(\bm x) < 0, \\
\displaystyle{\frac{\pi}{2 \alpha} \sin \left( \pi \frac{ q}{\alpha} \right)}, & 0 \leq q(\bm x) < \alpha \\
0, & q(\bm x) \geq \alpha,
\end{cases}
\end{equation} 
where $\alpha > 0$  is smoothing parameter.

\begin{figure}
\begin{center}
\includegraphics[width=0.8\textwidth]{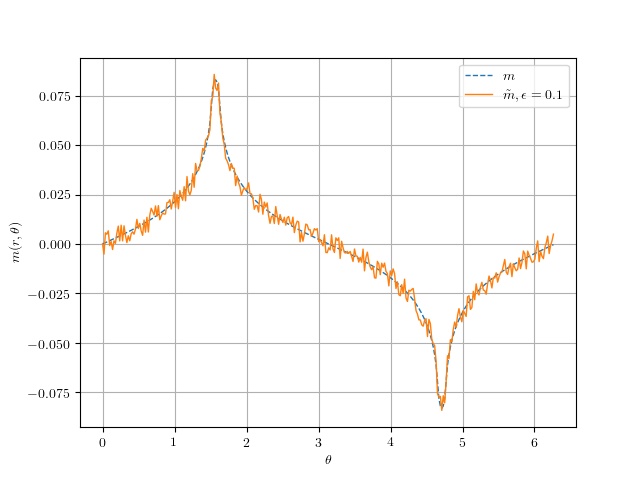}
\end{center}
\caption{Noiseless data $m_1(\bm x)$ (blue dashed line) and noisy data $\tilde{m}_1(\bm x)$ with $\epsilon=0.1$ (orange solid line)}
\label{fig:noisy_data}
\end{figure}

We add some uniform noise to the generated data.  More precisely, the noisy data $\tilde{m}_j(\bm x)$ is obtained by adding to $m_j(\bm x)$ a uniform  noise as follows:
\[
\tilde{m}_j(\bm x) = m_j(\bm x)  + \epsilon \|m_j(\bm x) \|_{L^2(\partial \Omega)} \frac{\theta_j(\bm x) }{\|\theta_j(\bm x) \|_{L^2(\partial \Omega)} }, \, j = 1, \dots, M, \, \bm x \in \partial \Omega,
\]
where $\epsilon$ is the noise level, $\theta_j (\bm x)$ is random numbers uniformly distributed on the interval $(-1.0, 1.0)$. Note that the noise is added only on the boundary $\partial \Omega$. The norms $ \|m_j(\bm x) \|_{L^2(\partial \Omega)}$ and $\|\theta_j(\bm x) \|_{L^2(\partial \Omega)}$ are used to scale $\theta_j (\bm x)$. Figure \ref{fig:noisy_data} shows noiseless and noisy data with $\epsilon=0.1$ on the perimeter of the computational domain $\partial \Omega$. Noiseless data is illustrated by blue dashed line, noisy data  by orange solid line. The data is generated for ellipse with 2 electrodes (1 measurement).

\begin{figure}
\begin{center}
\subfloat[\label{fig:f0}]{\includegraphics[width=0.33\textwidth]{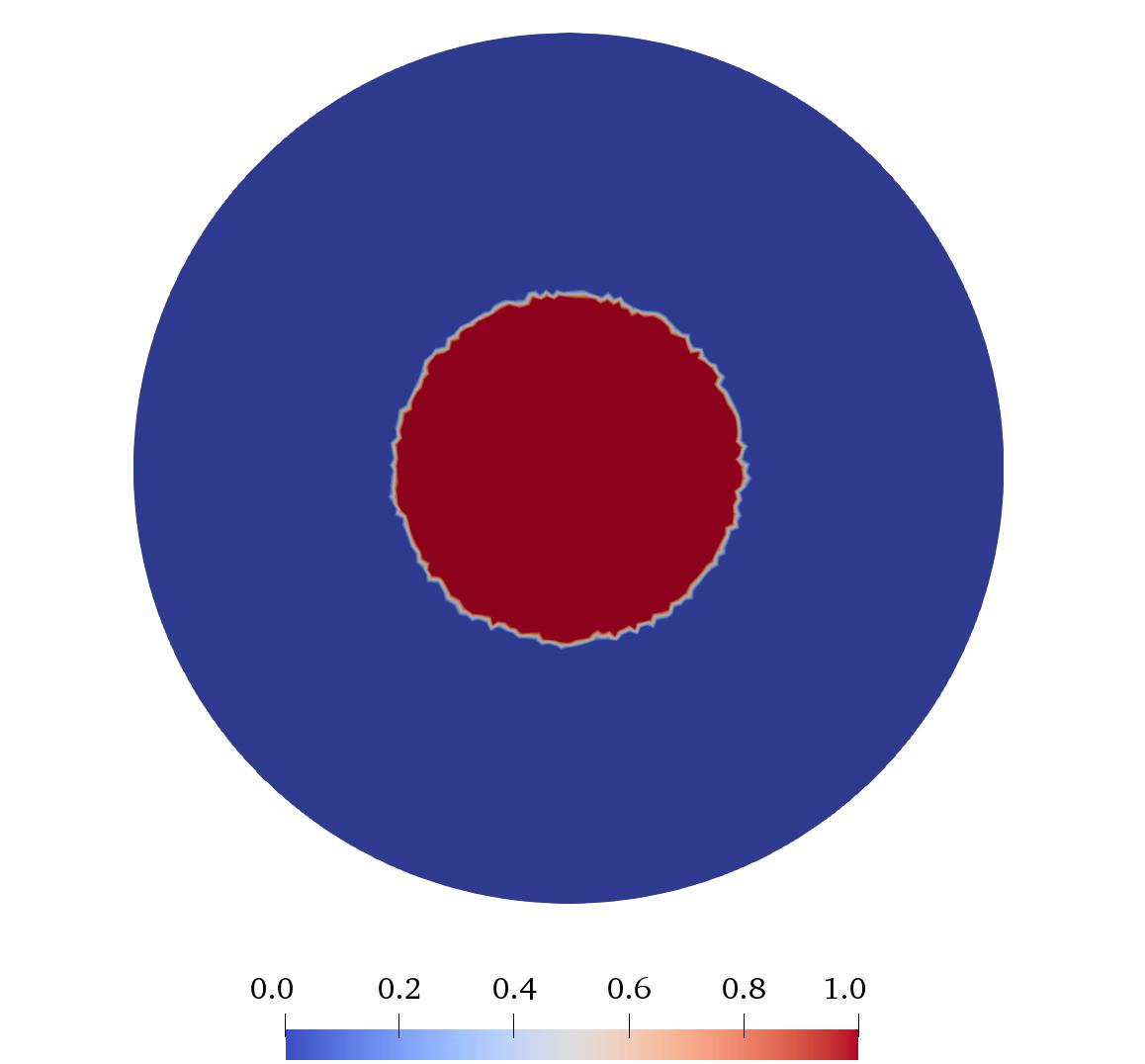}}
\subfloat[\label{fig:q0}]{\includegraphics[width=0.33\textwidth]{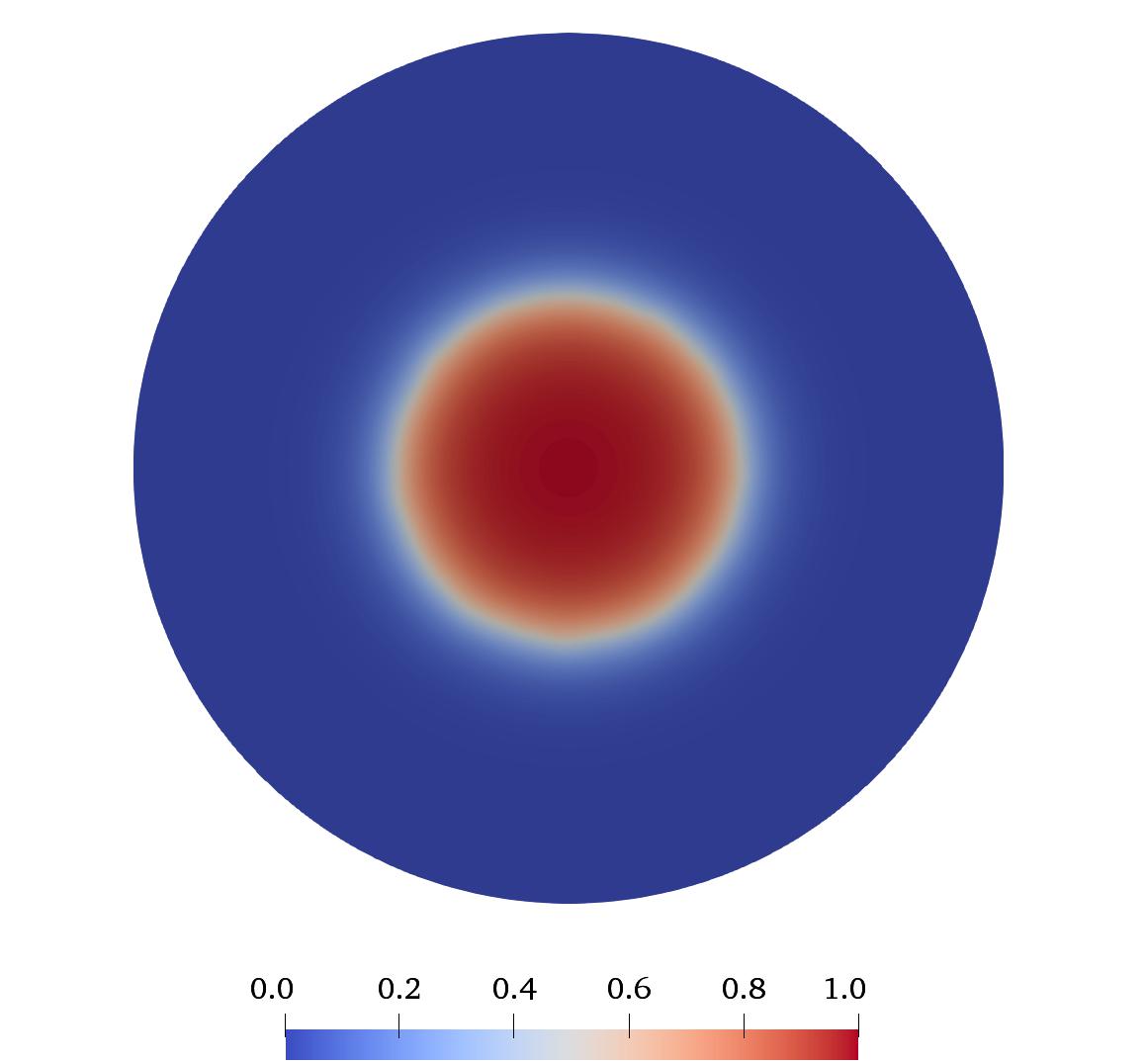}} %
\subfloat[\label{fig:c0}]{\includegraphics[width=0.33\textwidth]{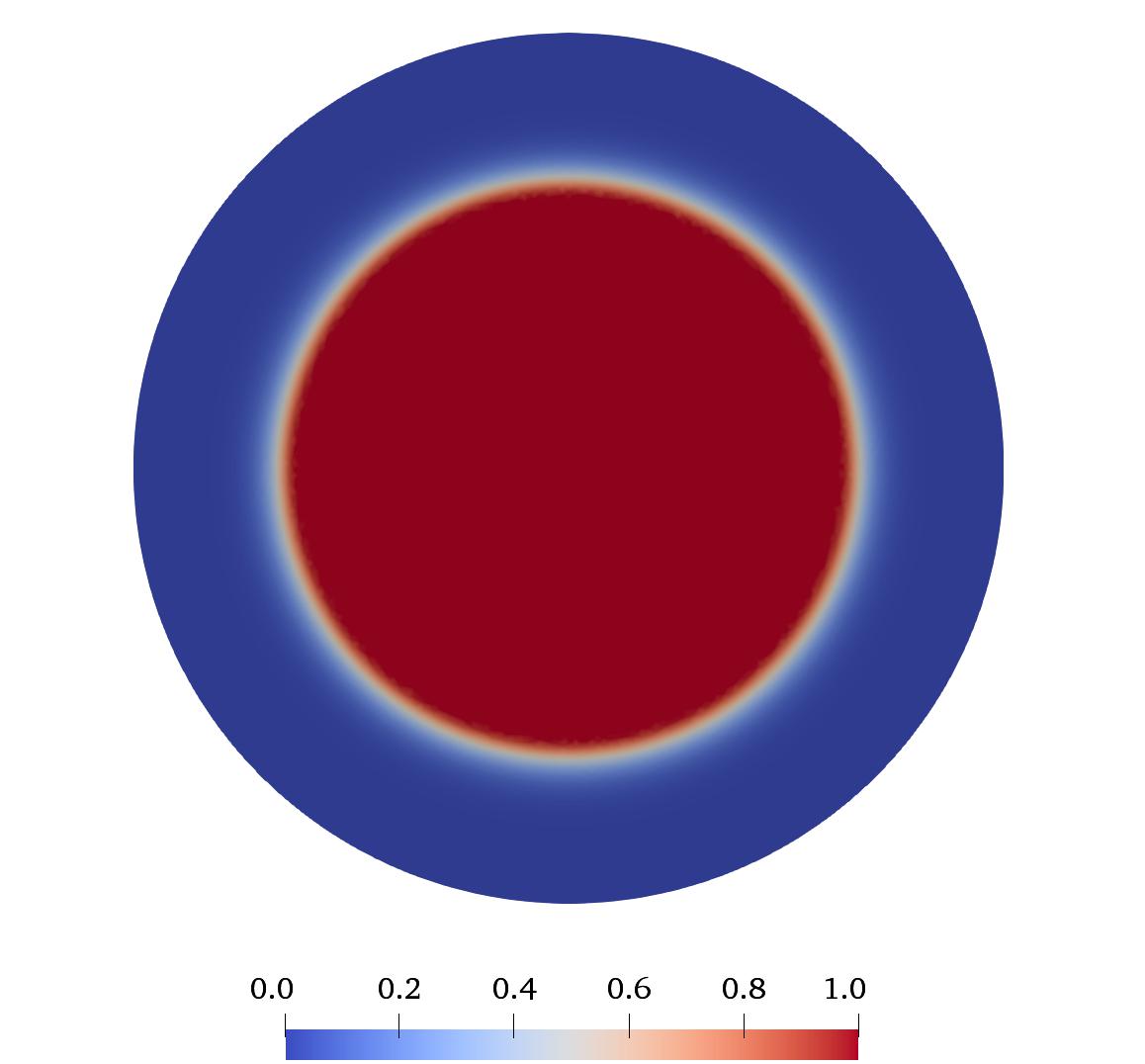}} 
\end{center}
\caption{ The initial guess for: (a) the control function $f^0(\bm x)$,  (b) the function $q^0(\bm x)$,  (c) the smoothed Heaviside function $H^0_{\alpha}(q(\bm x))$  }
\label{fig:initial_fq}
\end{figure}

The initial guess for the control function $f(\bm x)$ is chosen as follows
\begin{equation*}
f^0(\bm x) = 
\begin{cases}
1, & \bm x \in D^0 , \\
0, & \bm x \in \Omega \backslash D^0,
\end{cases}
\end{equation*}
where $D^0$ is some initial guess. It is taken as a circle  with the radius $r = 0.2$ and the center $\bm x_0 = (0.0, 0.0)$ of the domain $\Omega$. This initial control function $f^0(\bm x)$ is shown in Figure \ref{fig:f0}. 
In Figure \ref{fig:q0} we show the functions $q^0(\bm x)$ obtained by solving (\ref{eq:aux_var}) with $\gamma = 0.001$. 
The initial guess for the smoothed Heaviside function $H^0_{\alpha}(q(\bm x))$ (\ref{eq:heaviside_alpha}) with $\alpha = 0.01$ corresponding to $q^0(\bm x)$ in Figure \ref{fig:q0} is displayed in Figure \ref{fig:c0}. Note that $H^0_{\alpha}(\bm x)$ is shown in red (equal to 1)  for $q(\bm x) > 0.01$.

\begin{figure}
\begin{center}
\includegraphics[width=0.8\textwidth]{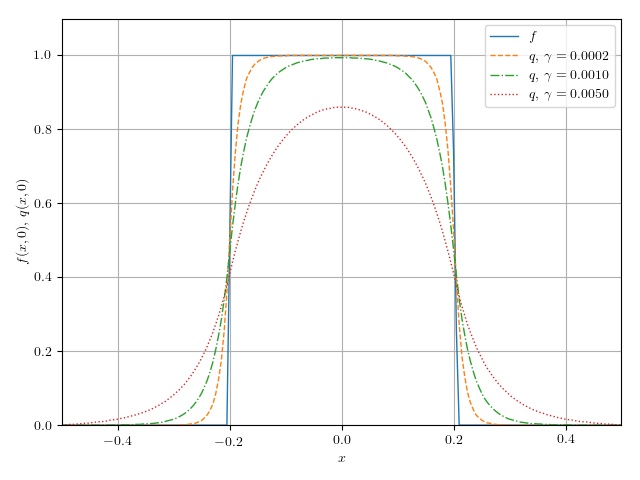}
\end{center}
\caption{Dependence of the level set function $q(\bm x)$ on the parameter $\gamma$}
\label{fig:fq_gammas}
\end{figure}

Figure \ref{fig:fq_gammas} shows the dependence of the level set function $q(\bm x)$ on the parameter $\gamma$. Here, the blue solid line plots  the values of $f(\bm x)$ over a line $(-0.5,0.0)$--$(0.5,0.0)$, the orange dashed  and green dash-dotted lines plot the values of $q(\bm x)$, corresponding to $\gamma=0.0005$, $\gamma=0.001$, and $\gamma=0.005$, respectively, over the same line. We can see that  as the value of $\gamma$ increases, the level set functions $q(\bm x)$ is more smoothed out. 

For minimization, we use the steepest descent method.  
The stopping criterion depends on the noise level $\epsilon$ as follows:
\[
\left \| \frac{dJ}{df}^k \right \|_{L^{\infty}(\Omega)} \leq \epsilon \beta,
\]
where $k$ is iteration number, $\beta$ is equal to $10^{-9}$.

\subsection{Reconstruction of ellipse}

First, we consider the ellipse object  (Figure \ref{fig:mesh_ellipse}). 
We want to examine the dependence of reconstruction results on the number of measurements $M$ and the parameter $\gamma$. 
For this, we reconstruct the piecewise constant coefficient from computationally simulated data with uniform noise with $0.01$ noise level. 
In Table \ref{tab:ellipse_Je} the first and second columns show the number of measurements $M$ and the parameter $\gamma$, respectively. In the next two columns we show the number of iterations until convergence and the corresponding value of the functional $J$. The last column represent the reconstruction error $\varepsilon$ calculated as follows:
\begin{equation}
\label{eq:err}
\varepsilon = \frac{\| \chi_D(\bm x) - H_{\alpha}(q(\bm x)) \|_{L^2(\Omega)}}{\|\chi_D(\bm x) \|_{L^2(\Omega)}},
\end{equation}
where $\chi_D (\bm x)$ is the characteristic function of ellipse (\ref{eq:chi}) and $H_{\alpha} (q(\bm x))$ is the reconstructed value of the smoothed Heaviside function (\ref{eq:heaviside_alpha}). 
As can be seen in Table \ref{tab:ellipse_Je}, the number of measurements $M$ is crucial for reconstruction accuracy. Using only one measurement $M=1$ is clearly not enough to accurately recover the desired coefficient and the number of iteration until convergence, in this case, is significantly larger. In this case, three measurements are sufficient for reconstruction, since using four measurements leads to slightly worse results. As for the parameter $\gamma$, we need to use $\gamma = 0.001$ for the best results. Using less than or greater than $0.001$ leads to increased error $\varepsilon$ and more iterations are needed for convergence of minimization.  Note that in these calculations we use $\alpha = 0.01$.

\begin{figure}
\begin{center}
\subfloat[\label{fig:j_ellipse_gammas}]{\includegraphics[width=0.5\textwidth]{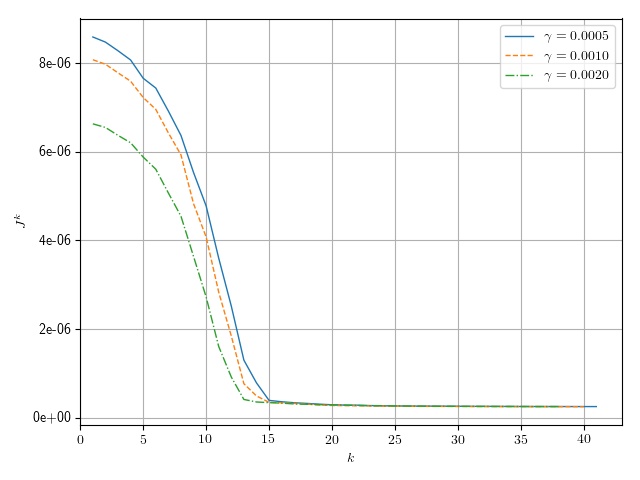}} 
\subfloat[\label{fig:e_ellipse_gammas}]{\includegraphics[width=0.5\textwidth]{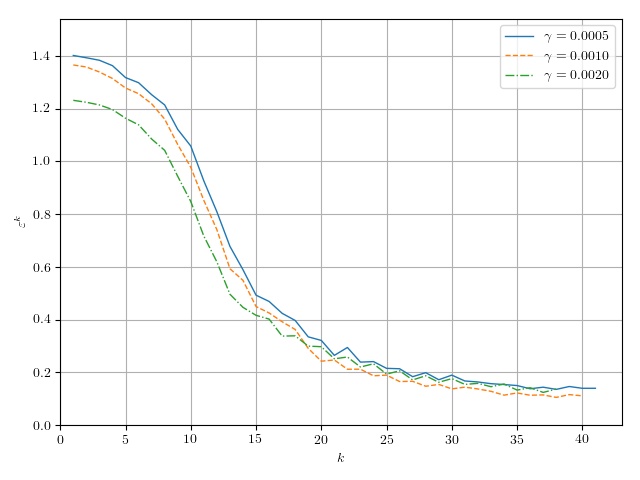}} 
\end{center}
\caption{The evolution of functional $J^k$ and error $\varepsilon^k$ with iterations $k$ for different values of $\gamma$: (a)  -- $J^k$,  (b)   -- $\varepsilon^k$}
\label{fig:j_err_ellipse_gammas}
\end{figure}

\begin{figure}
\begin{center}
\subfloat[\label{fig:j_ellipse_Ms}]{\includegraphics[width=0.5\textwidth]{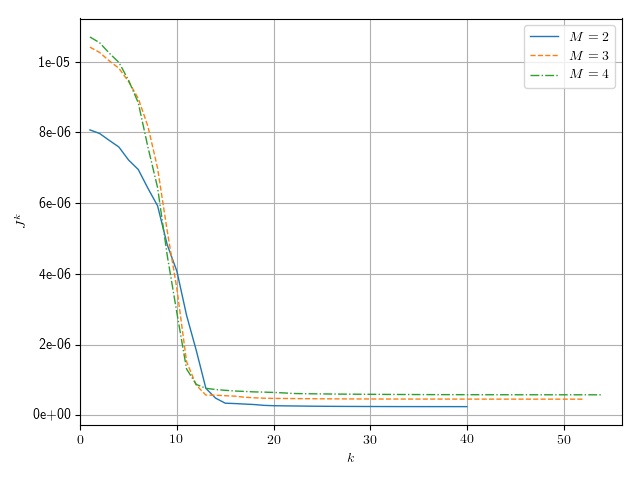}} 
\subfloat[\label{fig:e_ellipse_Ms}]{\includegraphics[width=0.5\textwidth]{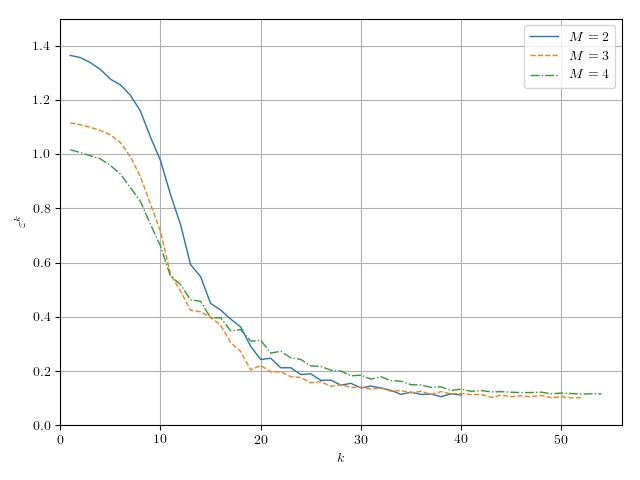}} 
\end{center}
\caption{The evolution of functional $J^k$ and error $\varepsilon^k$ with iterations $k$ for different values of $M$:  (a) -- $J^k$,  (b)  -- $\varepsilon^k$}
\label{fig:j_err_ellipse_Ms}
\end{figure}

In Figures \ref{fig:j_err_ellipse_gammas} and \ref{fig:j_err_ellipse_Ms}, we show the evolution of functional $J^k$ and error $\varepsilon^k$ with iterations $k$ for different values of $\gamma$ and $M$, respectively. 
In Figure \ref{fig:j_err_ellipse_gammas}, the blue solid, orange dashed, and green dash-dotted lines shows the values of $J^k$ (Figure \ref{fig:j_ellipse_gammas}) and $\varepsilon^k$   (Figure \ref{fig:e_ellipse_gammas}) for  $\gamma = 0.0005$,  $ 0.001$,  and $0.002$, respectively. 
Similarly, in Figure \ref{fig:j_err_ellipse_Ms}, the blue solid, orange dashed, and green dash-dotted lines shows the values of $J^k$ (Figure \ref{fig:j_ellipse_Ms}) and $\varepsilon^k$   (Figure \ref{fig:e_ellipse_Ms}) for  $M = 2$,  $3$,  and $4$, respectively. 
Reconstructions with different  $\gamma$ are performed with 3 measurements ($M = 3$).
While, $\gamma = 0.001$ is used for reconstructions with different $M$. 
In Figure \ref{fig:j_err_ellipse_gammas}, we see that the larger the value of $\gamma$, the smaller the initial values of the functional $J$ and the error $\varepsilon$. 
Although after a certain amount of iterations, the reconstruction error $\varepsilon$ for $\gamma=0.002$ becomes larger than the error $\varepsilon$ for $\gamma=0.001$ (see  Figure \ref{fig:e_ellipse_gammas}).  
This is due to the fact that the larger $\gamma$,  the smoother the function $q(\bm x)$, and this results in loss of detail in the smoothed Heaviside function $H_{\alpha}(\bm x)$. Similar situation  can be seen in Figure \ref{fig:j_err_ellipse_Ms}. When we use more measurements,  the initial error $\varepsilon$ increases. But again, in the process of minimization, the error $\varepsilon$ for $M=4$ becomes larger than the error $\varepsilon$ for $M=3$. We can also see that the larger the number of measurements $M$, the minimization process takes more iterations to converge. 

\begin{figure}
\begin{center}
\subfloat[\label{fig:ellipse_iter1}]{\includegraphics[width=0.33\textwidth]{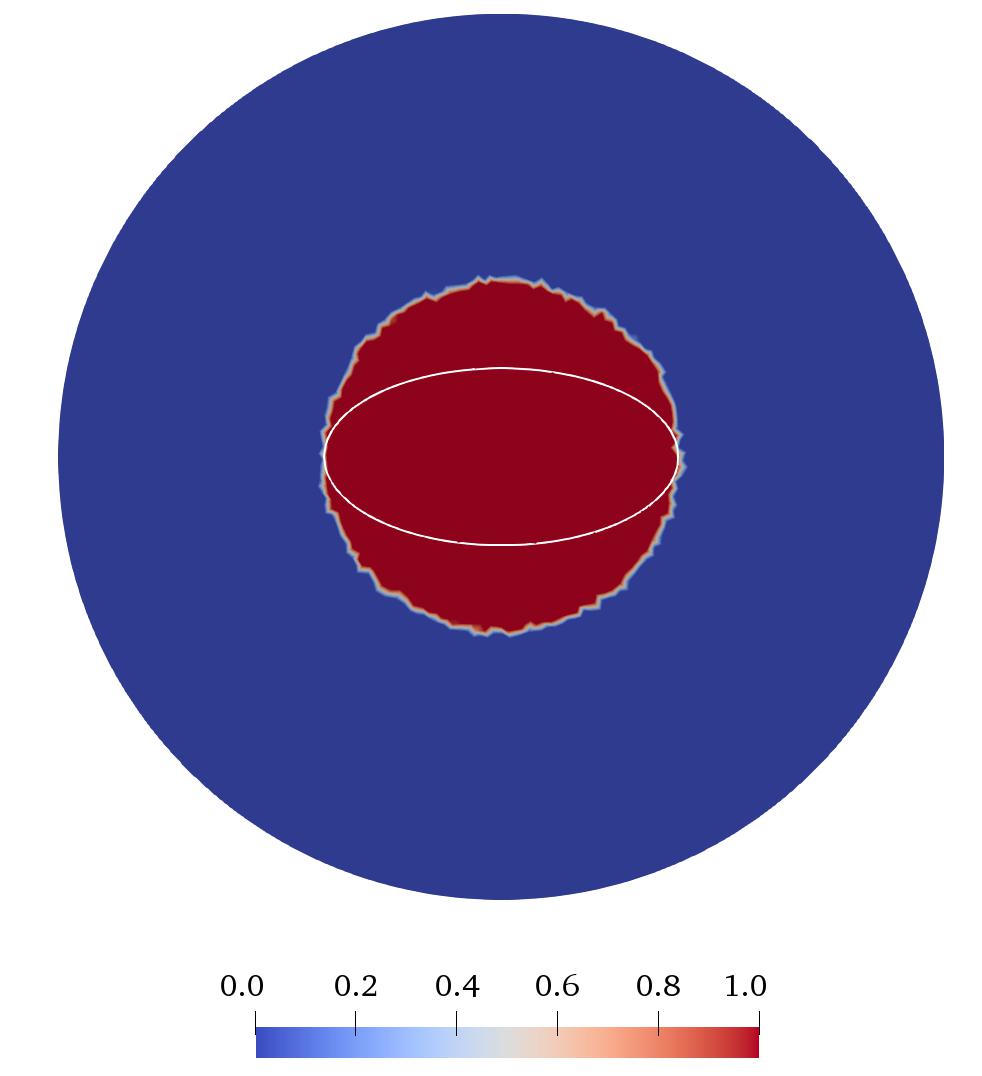}} %
\subfloat[\label{fig:ellipse_iter10}]{\includegraphics[width=0.33\textwidth]{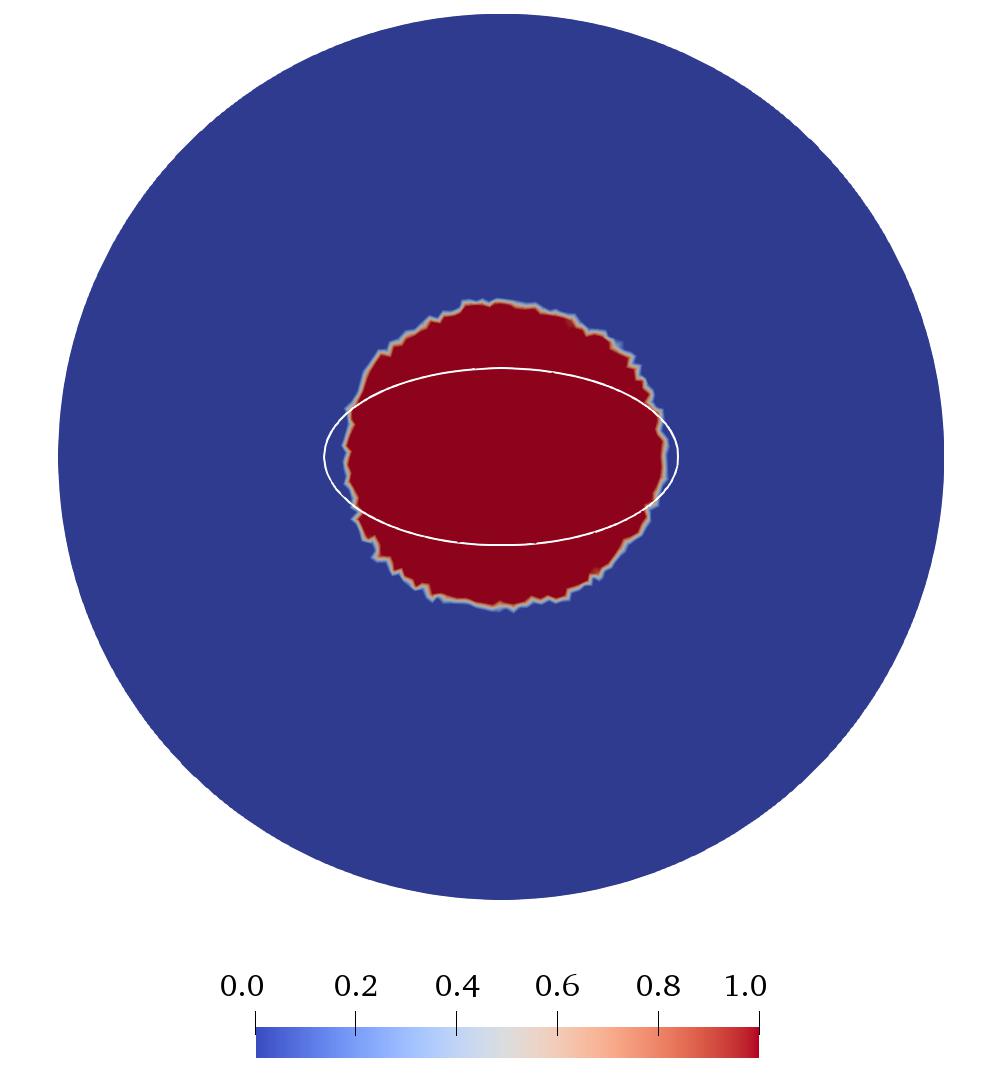}}  
\subfloat[\label{fig:ellipse_iter15}]{\includegraphics[width=0.33\textwidth]{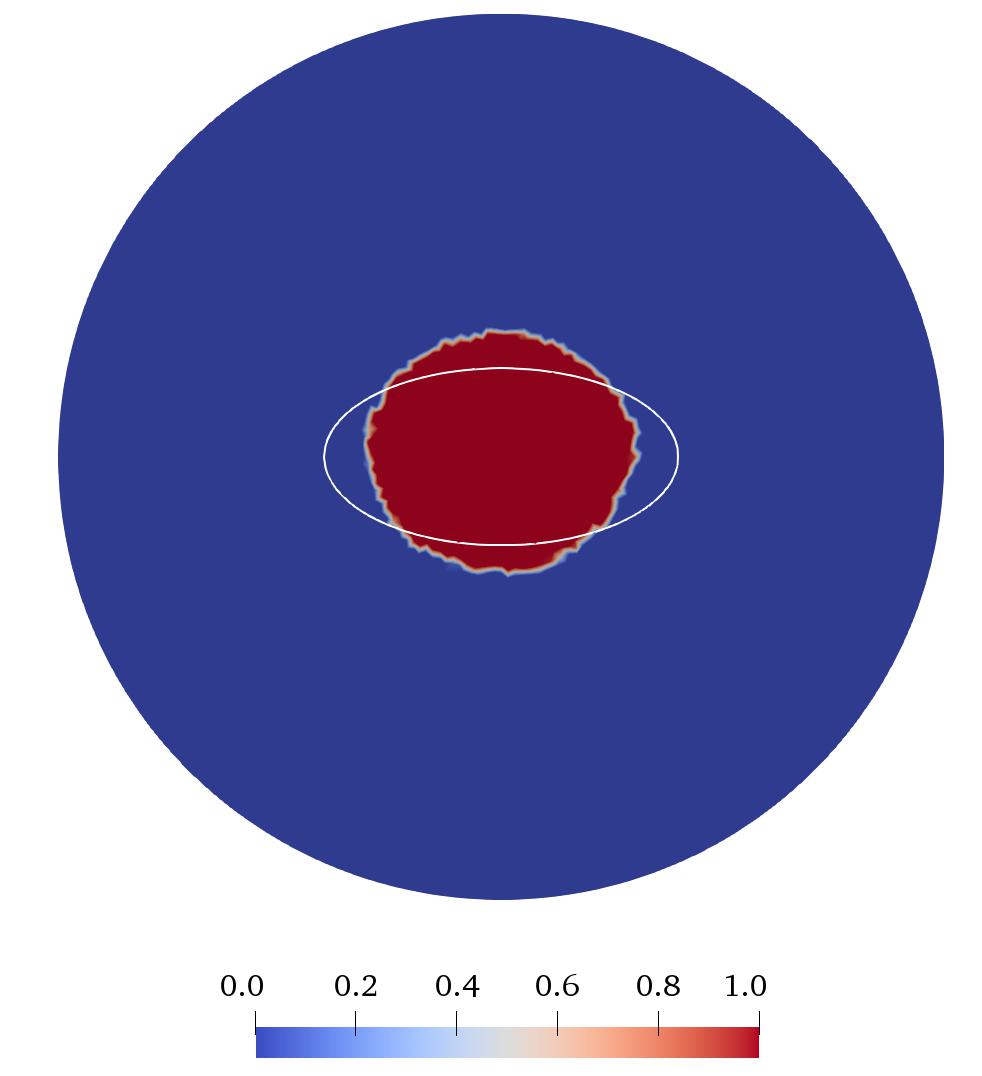}} \\
\subfloat[\label{fig:ellipse_iter20}]{\includegraphics[width=0.33\textwidth]{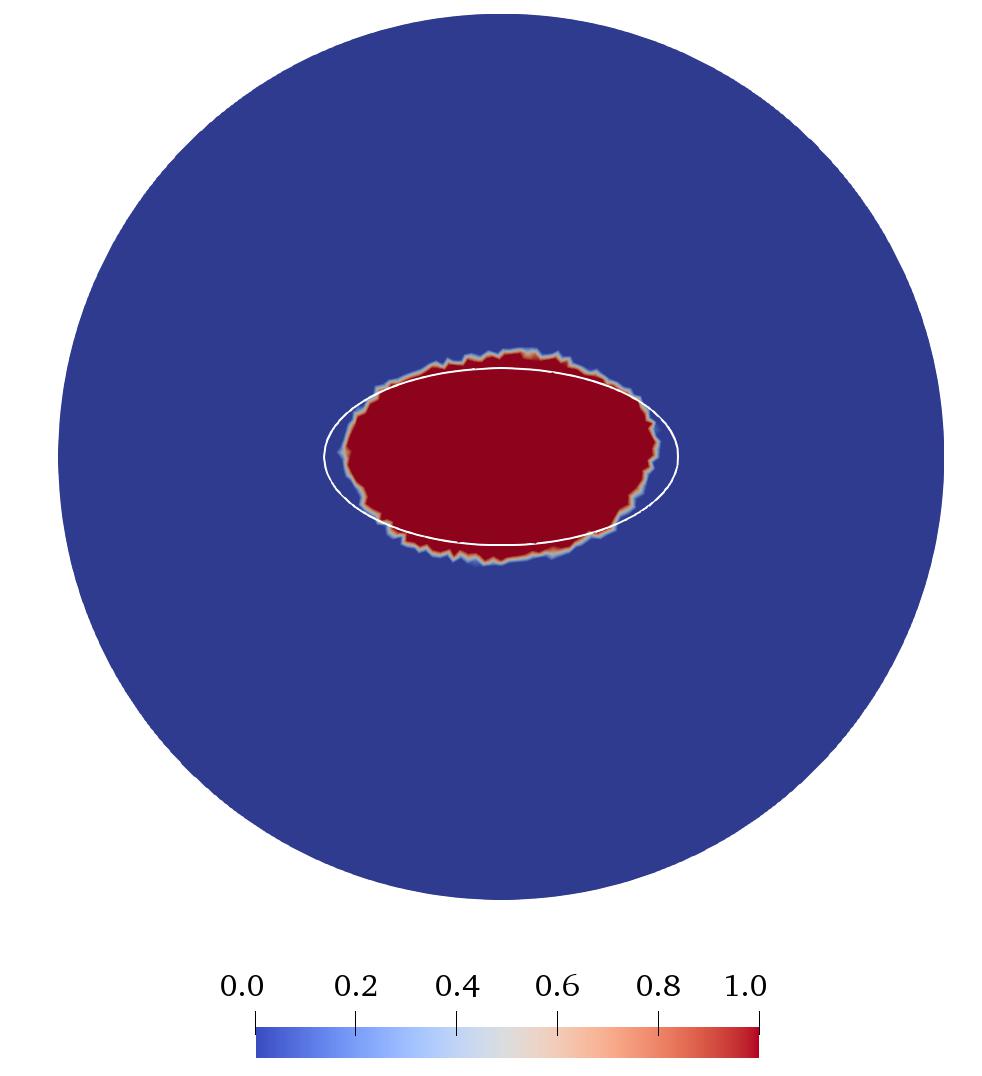}} 
\subfloat[\label{fig:ellipse_iter30}]{\includegraphics[width=0.33\textwidth]{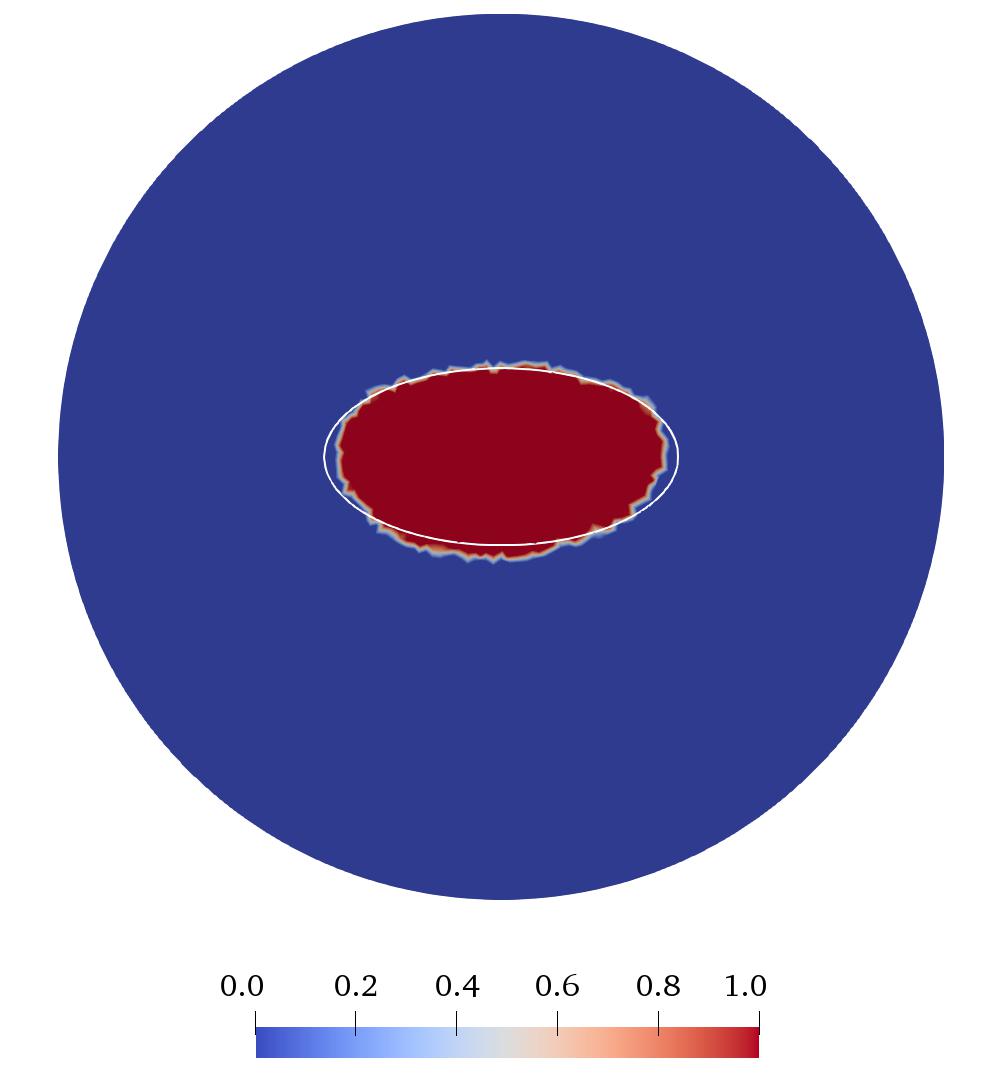}} 
\subfloat[\label{fig:ellipse_iter50}]{\includegraphics[width=0.33\textwidth]{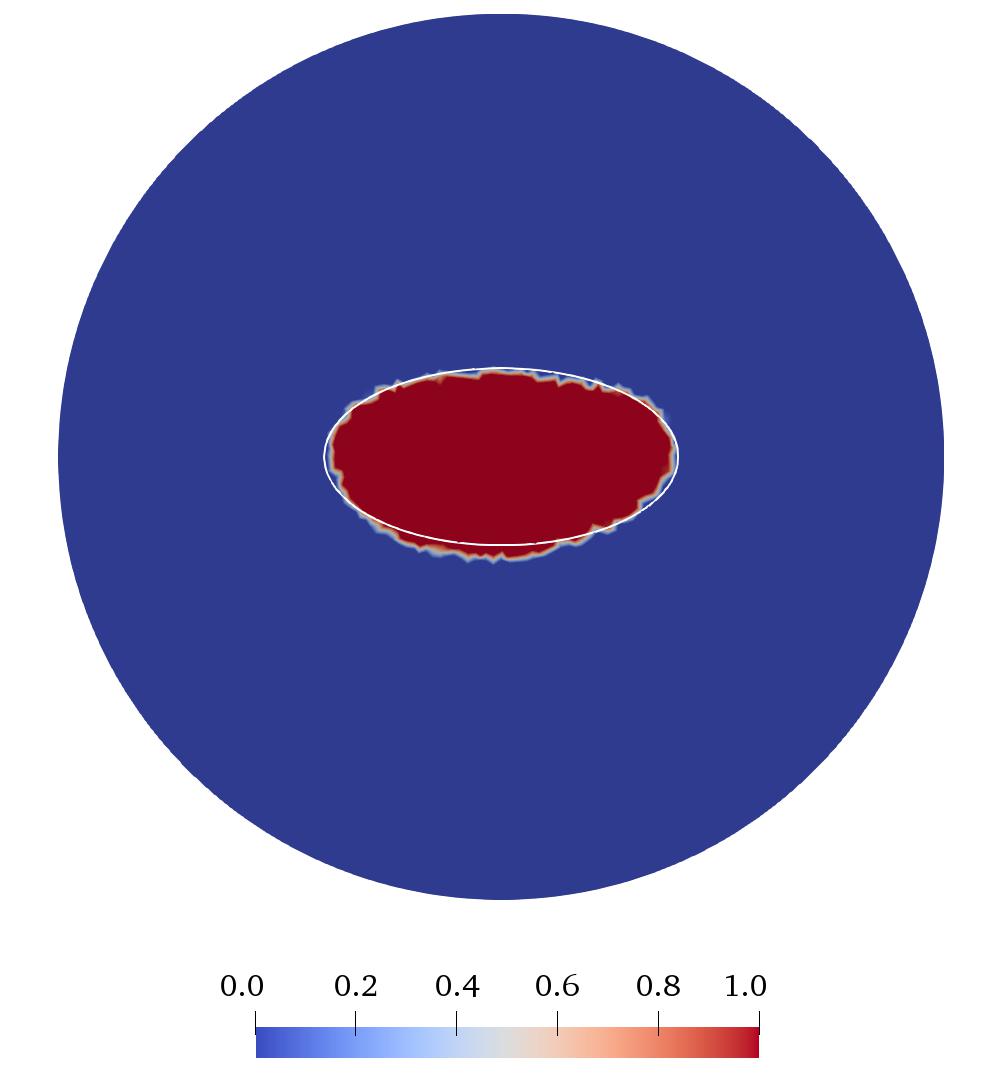}} 
\end{center}
\caption{Reconstructed smoothed Heaviside function $H_{\alpha}(q(\bm x))$  after: (a) --  1, (b) -- 10,  (c) -- 15,  (d) -- 20,  (e) -- 30, and (f) -- 50  iterations }
\label{fig:ellipse_iters}
\end{figure}

Figure \ref{fig:ellipse_iters} shows the reconstructed smoothed Heaviside functions $H_{\alpha}(q(\bm x))$ after 1, 10, 15, 20, 30 and 50 iterations. These results are obtained from 3 measurements with noise $\epsilon=0.01$ using $\gamma = 0.001$ and $\alpha = 0.01$. The shape of the true object (ellipse) is outlined with the solid white line. Note that in Figures  \ref{fig:j_err_ellipse_gammas} and \ref{fig:j_err_ellipse_Ms} the reconstruction results correspond to the orange dashed line. We can see that during the first 10 iterations, the recovered coefficient does not change much, and starting from 20 iterations, we already recovered the ellipse shape. Subsequent iterations are needed to increase the quality of recovery. 

\begin{table}
\begin{center}
\begin{tabular}{c|c|c|c|c}
\hline
$M$ & $\gamma$ & Iterations & $J\cdot 10^{-7}$ & $\varepsilon$ \\ \hline
	& 0.0005	& 308	& 1.203	  & 0.225	 \\
1	& 0.0010	& 328	& 1.204	  & 0.211	 \\
	& 0.0020	& 351	& 1.205	  & 0.232	  \\ 
\hline
 	& 0.0005	& 43	& 2.478 & 0.138	  \\
2	& 0.0010	& 42	& 2.465 & 0.114	  \\
	& 0.0020	& 40	& 2.484  & 0.127	  \\
\hline
	& 0.0005	& 52	& 4.582	  & 0.102	  \\
3	& 0.0010	& 54	& 4.580     & 0.101	  \\
	& 0.0020	& 54	& 4.603	  & 0.120	  \\ 
\hline
	& 0.0005	& 46	& 5.839	  & 0.156	  \\
4	& 0.0010	& 56	& 5.808	  & 0.116 \\
	& 0.0020	& 63	& 5.826	  & 0.125	 \\
\hline
\end{tabular}
\end{center}
\caption{Dependence of $J$, $\varepsilon$, and number of iterations on the number of measurements $M$ and the parameter $\gamma$ for reconstruction of ellipse}
\label{tab:ellipse_Je}
\end{table}

\subsection{Reconstruction of circles} 

Next, we consider two circles  in Figure \ref{fig:mesh_circles}. 
We perform the same experiments as with the ellipse using computationally simulated data with uniform noise with noise level. 

In Table \ref{tab:circles_Je} we list the number of measurements $M$, the parameter $\gamma$, the number of iterations until convergence, the final value of the functional $J$, the reconstruction error $\varepsilon$ (\ref{eq:err}). 
We can confirm that the number of measurements $M$ is essential for good accuracy.
In this case, we need to use more than two measurements. 
When we use only one measurements $M=1$,  we get poor quality results with big errors.  Also,  the algorithm did not converged in 1000 iterations and was terminated.  
We can see that $\gamma = 0.006$ is the optimal value for better results and using less than or greater than it leads to increased error $\varepsilon$. 
Note that when we increase the value of $\gamma$, more iterations are needed for convergence of the algorithm.
Here, in all calculations we use $\alpha = 0.01$. 

\begin{figure}
\begin{center}
\subfloat[\label{fig:j_circles_gammas}]{\includegraphics[width=0.5\textwidth]{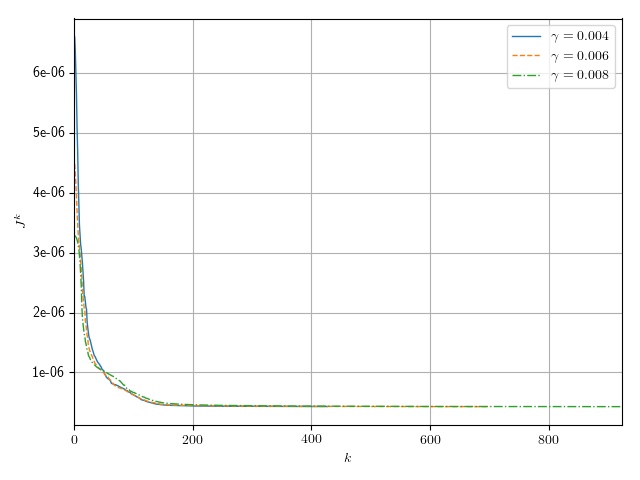}} 
\subfloat[\label{fig:e_circles_gammas}]{\includegraphics[width=0.5\textwidth]{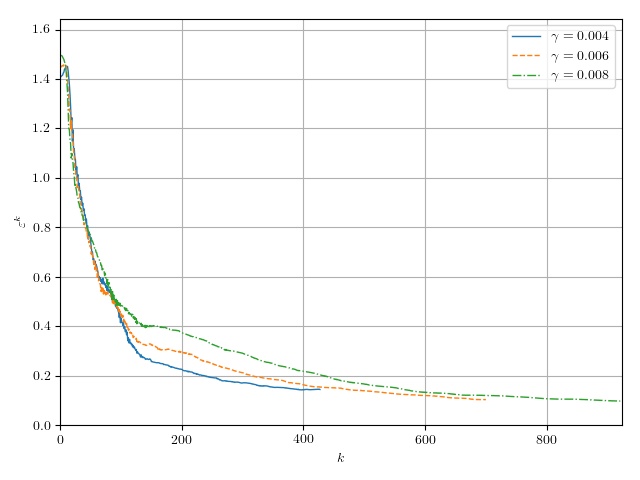}} 
\end{center}
\caption{The evolution of functional $J^k$ and error $\varepsilon^k$ with iterations $k$ for different values of $\gamma$: (a) -- $J^k$, (b)   -- $\varepsilon^k$}
\label{fig:j_err_circles_gammas}
\end{figure}

\begin{figure}
\begin{center}
\subfloat[\label{fig:j_circles_Ms}]{\includegraphics[width=0.5\textwidth]{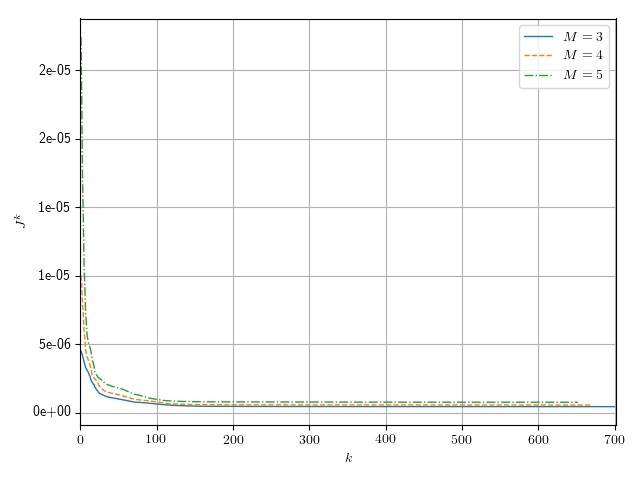}} 
\subfloat[\label{fig:e_circles_Ms}]{\includegraphics[width=0.5\textwidth]{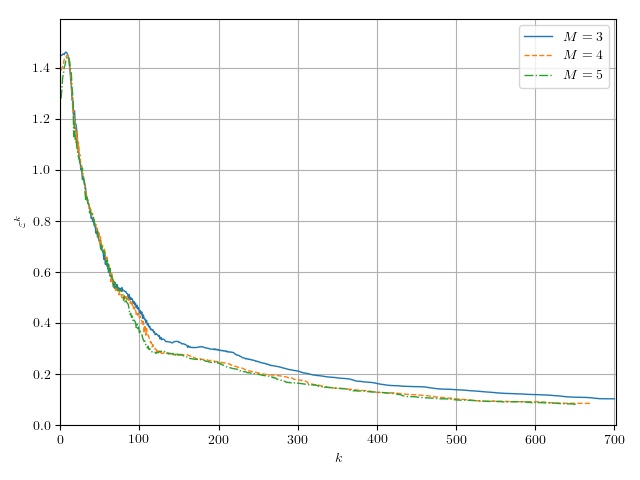}} 
\end{center}
\caption{The evolution of functional $J^k$ and error $\varepsilon^k$ with iterations $k$ for different values of $M$: (a) -- $J^k$,  (b)  -- $\varepsilon^k$}
\label{fig:j_err_circles_Ms}
\end{figure}

Similar to Figures \ref{fig:j_err_ellipse_gammas} and \ref{fig:j_err_ellipse_Ms} for the reconstruction of ellipse,
Figures \ref{fig:j_err_circles_gammas} and \ref{fig:j_err_circles_Ms} illustrate the evolution of functional $J^k$ and error $\varepsilon^k$ with iterations $k$ for different values of $\gamma$ and $M$, respectively. 
In Figure \ref{fig:j_err_circles_gammas}, 
the blue solid, orange dashed, and green dash-dotted lines shows 
the values of $J^k$ (Figure \ref{fig:j_circles_gammas}) and $\varepsilon^k$   (Figure \ref{fig:e_circles_gammas}) for  $\gamma = 0.004$,  $ 0.006$,  and $0.008$, respectively. 
Similarly, in Figure \ref{fig:j_err_circles_Ms}, the blue solid, orange dashed, and green dash-dotted lines shows the values of $J^k$ (Figure \ref{fig:j_circles_Ms}) and $\varepsilon^k$   (Figure \ref{fig:e_circles_Ms}) for  $M = 3$,  $4$,  and $5$, respectively. 
Reconstructions with different  $\gamma$ are performed with $M = 3$ and reconstructions with different $M$ are performed with $\gamma = 0.006$. 

\begin{figure}
\begin{center}
\subfloat[\label{fig:circles_iter5}]{\includegraphics[width=0.33\textwidth]{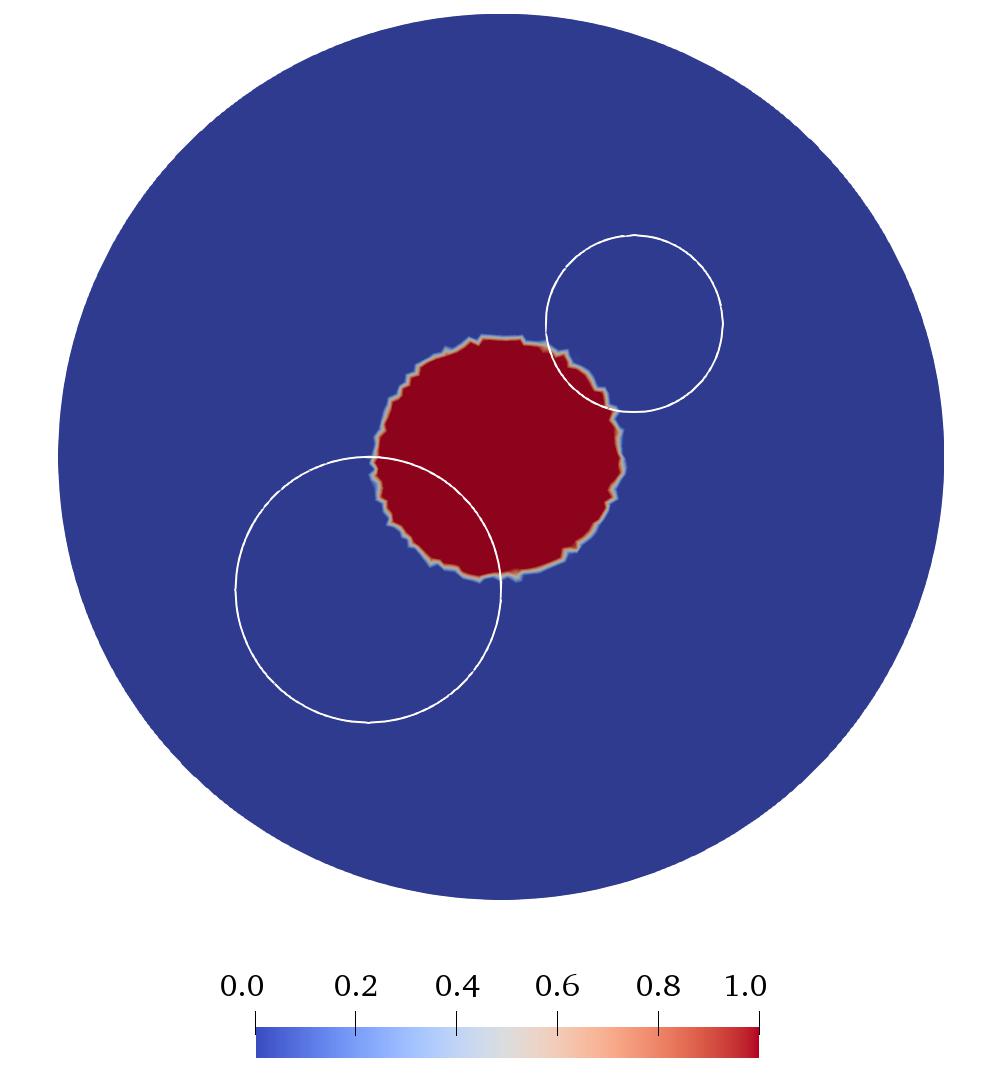}} %
\subfloat[\label{fig:circles_iter25}]{\includegraphics[width=0.33\textwidth]{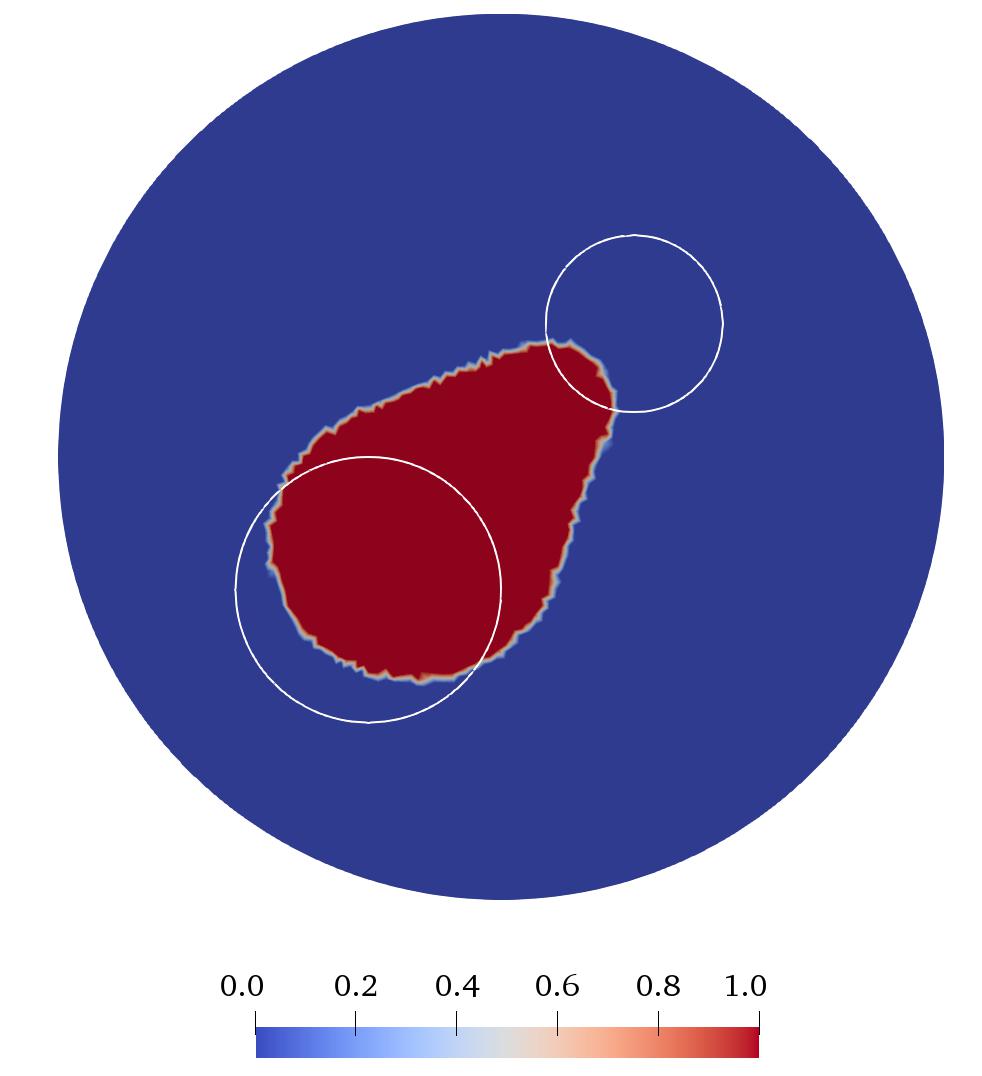}}  
\subfloat[\label{fig:circles_iter50}]{\includegraphics[width=0.33\textwidth]{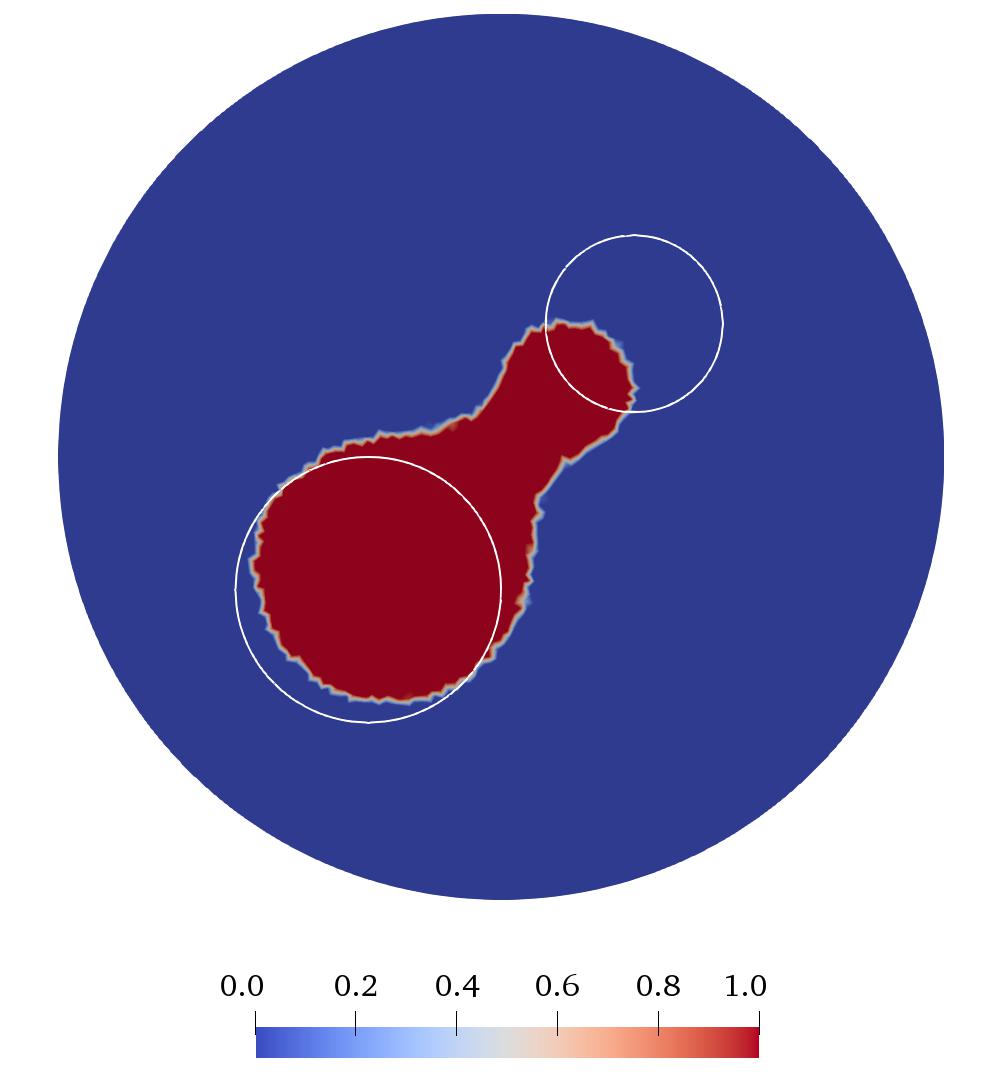}} \\
\subfloat[\label{fig:circles_iter150}]{\includegraphics[width=0.33\textwidth]{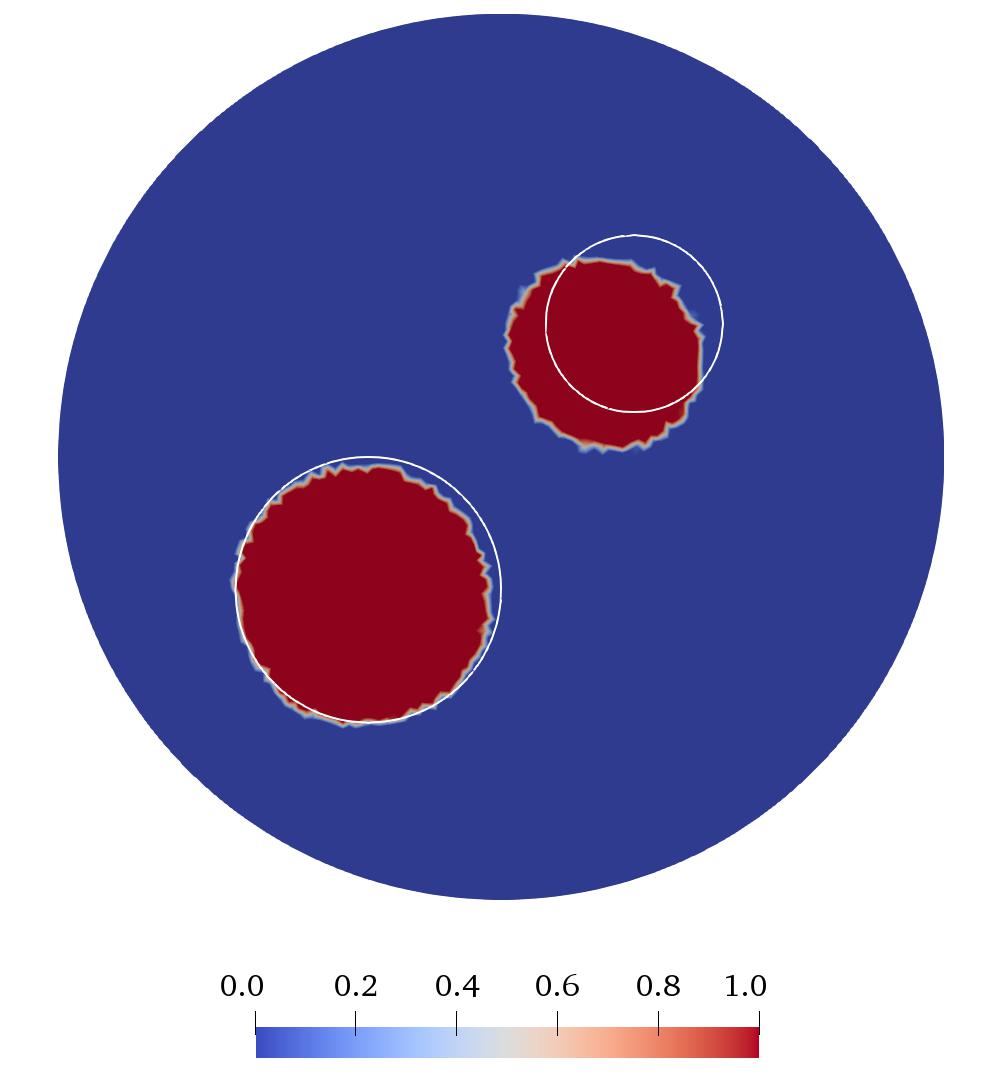}} 
\subfloat[\label{fig:circles_iter200}]{\includegraphics[width=0.33\textwidth]{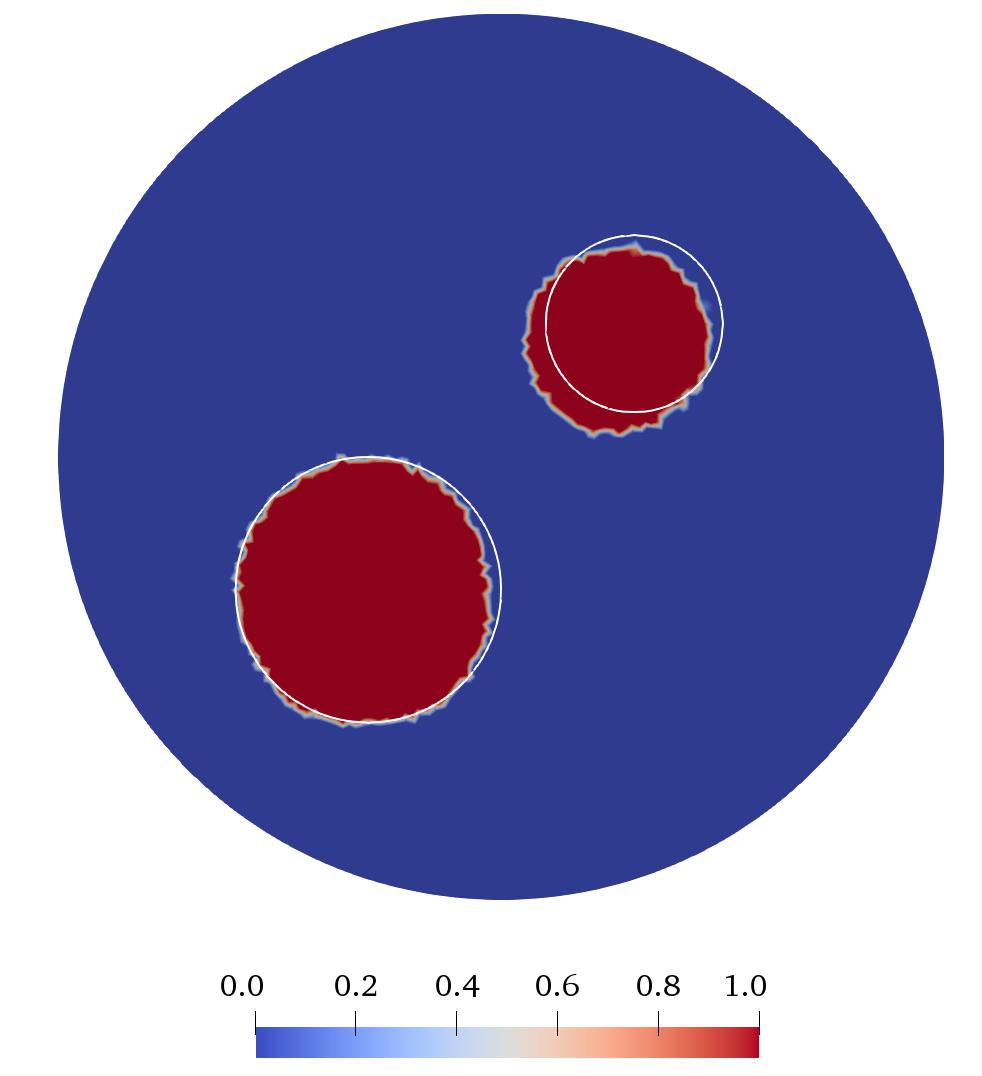}} 
\subfloat[\label{fig:circles_iter650}]{\includegraphics[width=0.33\textwidth]{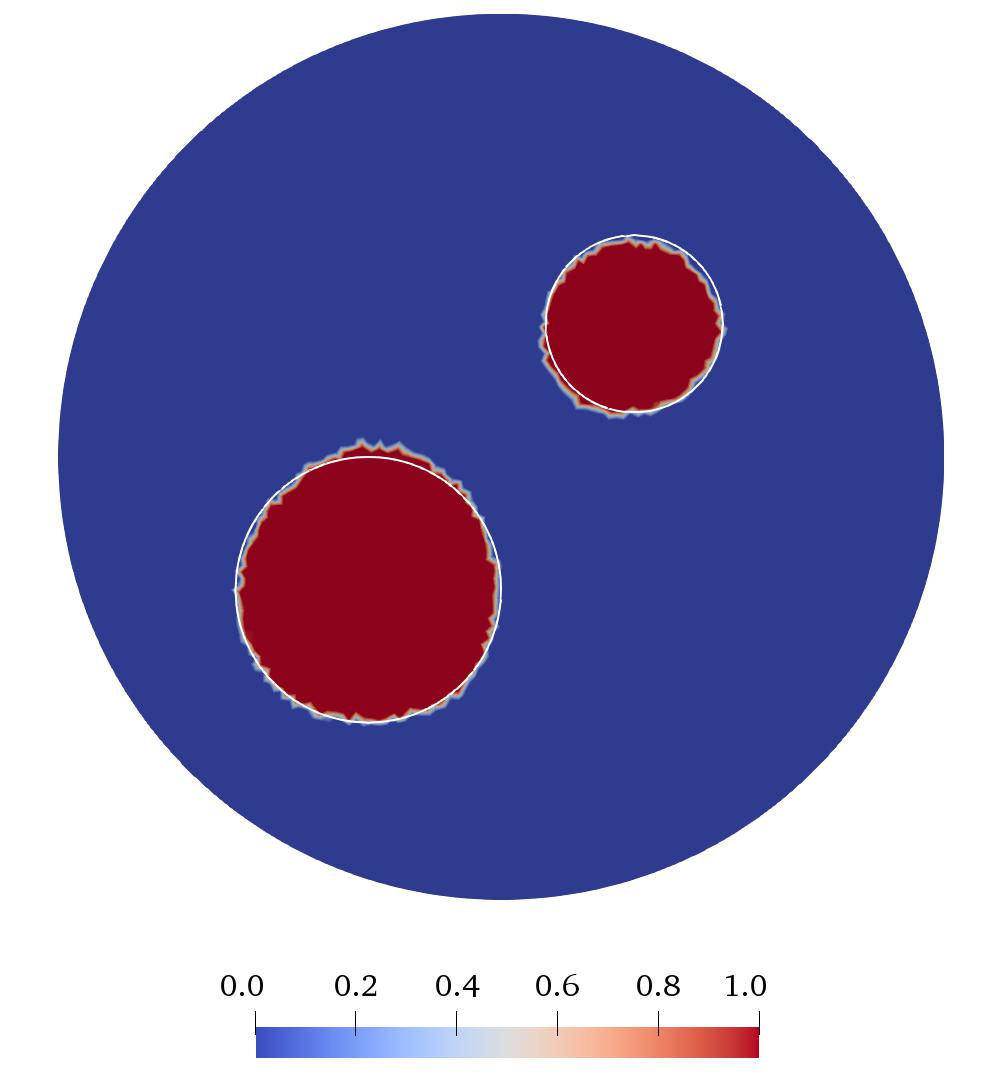}} 
\end{center}
\caption{Reconstructed smoothed Heaviside function $H_{\alpha}(q(\bm x))$  after: (a)  -- 5,  (b) -- 25, (c)  -- 50,  (d) -- 150,  (e) -- 300, and (f) -- 650  iterations }
\label{fig:circles_iters}
\end{figure}

Figure \ref{fig:circles_iters} shows the reconstructed smoothed Heaviside functions $H_{\alpha}(q(\bm x))$ after 5, 25, 50, 150, 300 and 650 iterations. These images are obtained from 5 measurements with noise $\epsilon=0.01$ using $\gamma = 0.006$ and $\alpha = 0.01$. The shape of true objects is outlined with the solid white line. We can see that after 150 iterations, the circles are almost recovered.  The location and shape of the bigger circle are already reconstructed. Further iterations are required to recover the true location of the second circle.

\begin{table}
\begin{center}
\begin{tabular}{c|c|c|c|c}
\hline
$M$ & $\gamma$ & Iterations & $J\cdot 10^{-7}$ & $\varepsilon$ \\ \hline
	& 0.004	& 1000	& 1.154 & 	  0.853	\\
1	& 0.006	& 1000	& 1.385 &  1.543 \\
	& 0.008	& 1000	& 1.702 &   1.645	\\
\hline
	& 0.004	& 92	    & 6.119 & 	  0.883	\\
2	& 0.006	& 112	& 6.113 & 	  0.863	 \\
	& 0.008	& 123	& 6.121 & 	  0.837	\\
\hline
	& 0.004	& 429	& 4.324	& 0.145	\\
3	& 0.006	& 702	& 4.307	& 0.103	\\
	& 0.008	& 923	& 4.306	& 0.097	\\
\hline
	& 0.004	& 300	& 5.546	& 0.149	\\
4	& 0.006	& 673	& 5.473	& 0.085	\\
	& 0.008	& 938	& 5.471		& 0.086	\\
\hline
	& 0.004	& 325	& 7.537		& 0.106	\\
5	& 0.006	& 654	& 7.503	& 0.082	\\
	& 0.008	& 956	& 7.496	& 0.084	\\
\hline
\end{tabular}
\end{center}
\caption{Dependence of $J$, $\varepsilon$, and number of iterations on the number of measurements $M$ and the parameter $\gamma$ for reconstruction of circles}
\label{tab:circles_Je}
\end{table}

\begin{figure}
\begin{center}
\subfloat[\label{fig:circles_err2}]{\includegraphics[width=0.33\textwidth]{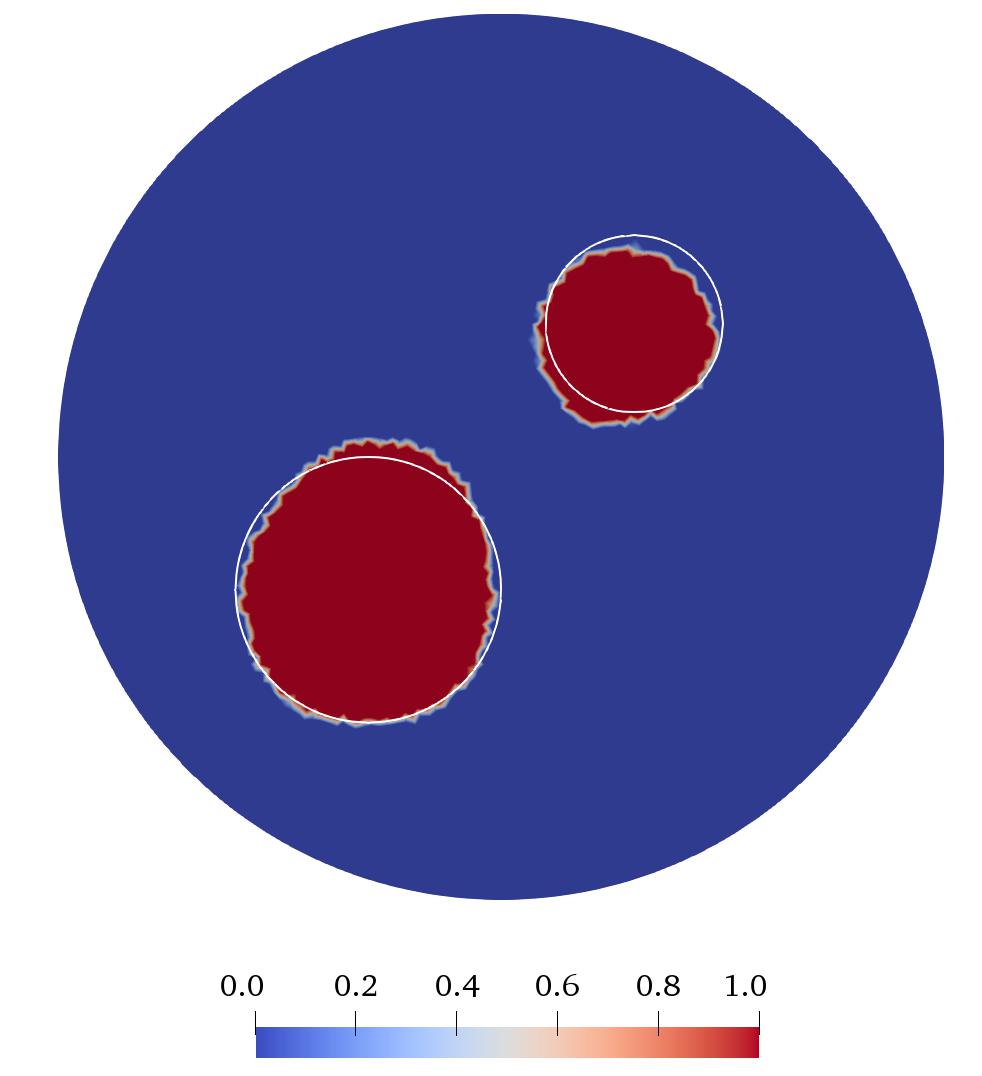}} %
\subfloat[\label{fig:circles_err3}]{\includegraphics[width=0.33\textwidth]{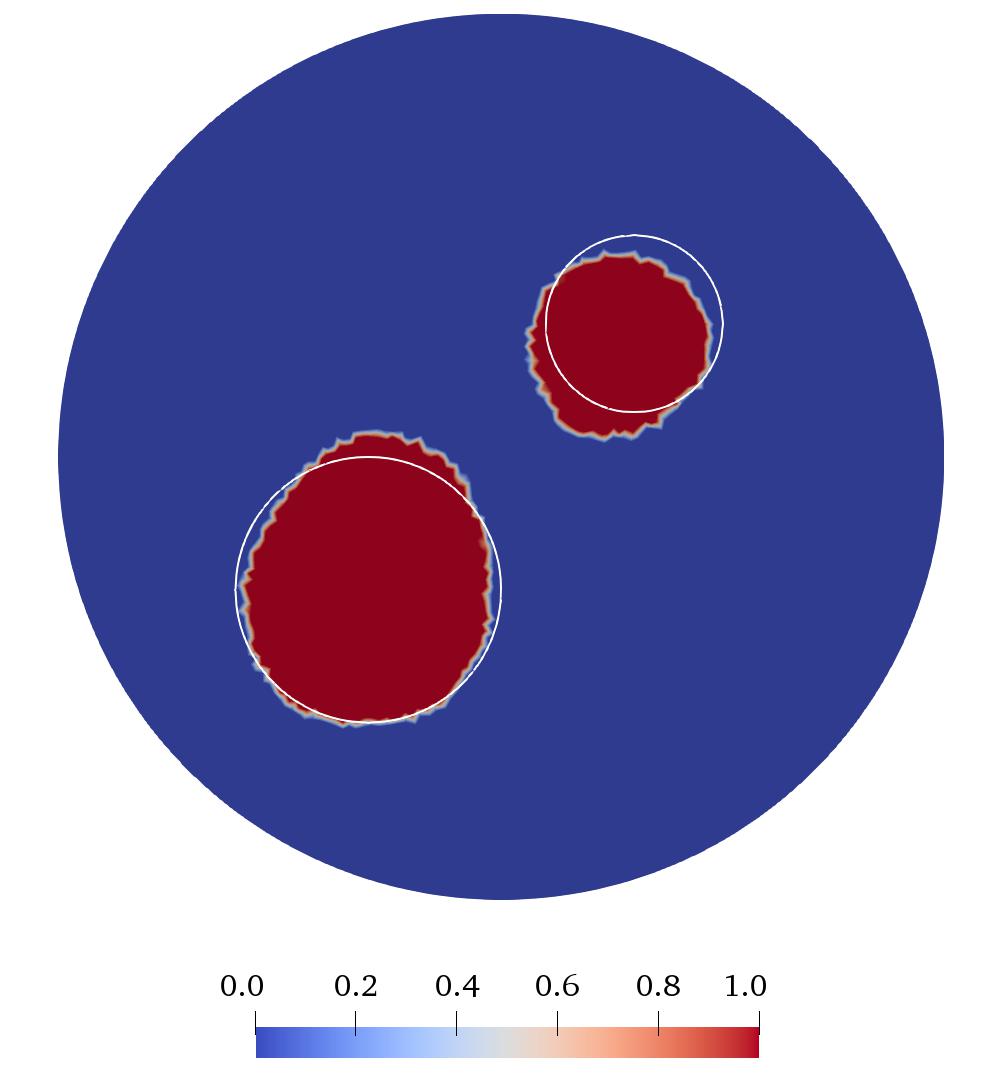}}  
\subfloat[\label{fig:circles_err4}]{\includegraphics[width=0.33\textwidth]{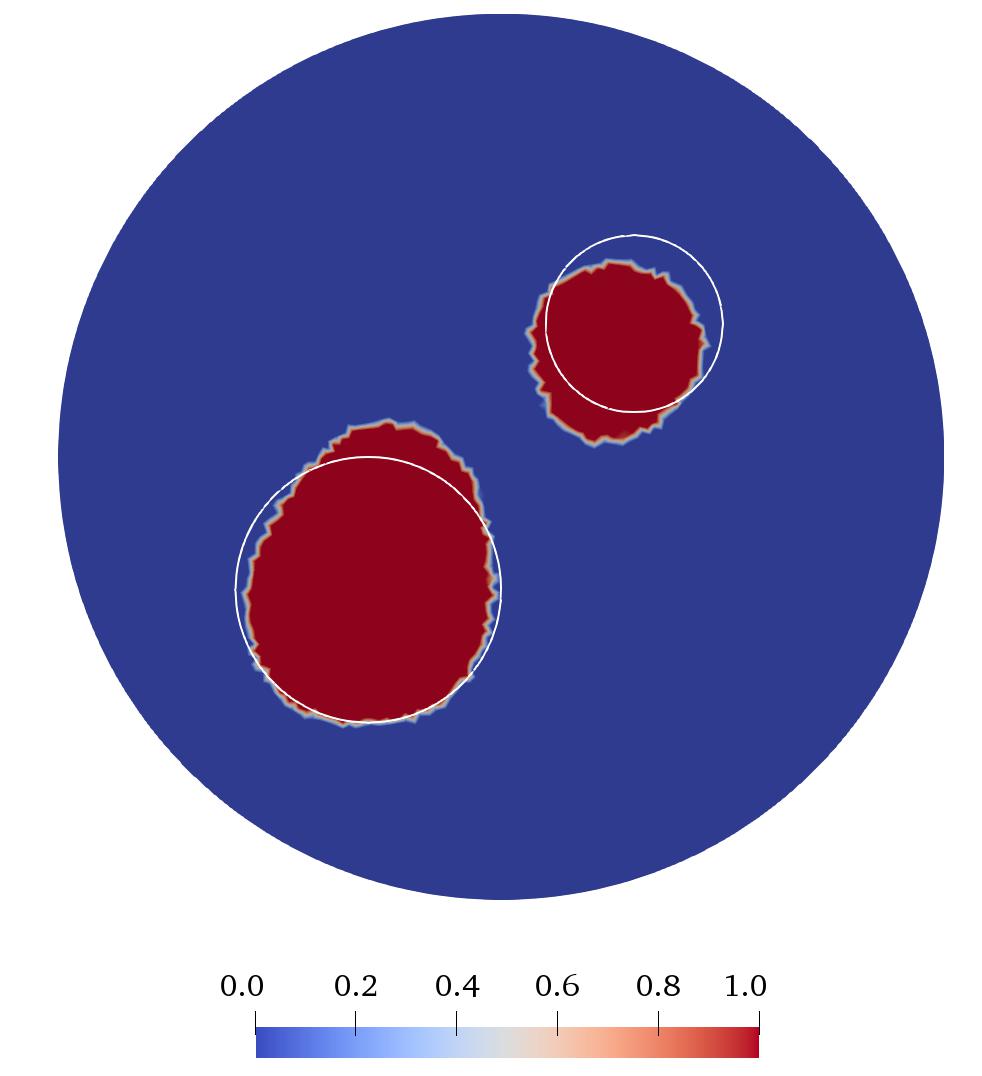}} 
\end{center}
\caption{Reconstructed smoothed Heaviside function $H_{\alpha}(q(\bm x))$  from data with noise: (a)  $\epsilon=0.02$,  (b) $\epsilon=0.03$,  (c) $\epsilon=0.04$, }
\label{fig:circles_err}
\end{figure}

Finally, we want to see how our algorithm can handle the different values of noise. We consider the two circles.
In Table \ref{tab:circles_err} the first and second columns show the number of measurements $M$ and the parameter $\gamma$, respectively. In the next four columns we show the reconstruction error $\varepsilon$ (\ref{eq:err}) for different values of noise level $\epsilon = 0.01$,  $0.02$, $0.03$ and $0.04$. 
As can be seen in Table \ref{tab:circles_err}, we need to increase the number of measurements $M$ to handle big noise in our data. 
Also, we observe that when $ M $ is less than 8, the effect of the value of $ \gamma $ 
That means that the effect of $M$ is much more significant for better reconstruction. 
When $M$ is equal or greater than 8, $\gamma = 0.006$ gives better accuracy. 
Figure \ref{fig:circles_err} illustrates the reconstructed smoothed Heaviside functions $H_{\alpha}(q(\bm x))$ from data with noise $\epsilon = 0.02$, $0.03$ and $0.04$. These results are obtained from 10 measurements using $\gamma = 0.006$.

\begin{table}
\begin{center}
\begin{tabular}{c|c|c|c|c|c}
\hline
$M$ & $\gamma$  & $\epsilon=0.01$ & $\epsilon=0.02$ & $\epsilon=0.03$ & $\epsilon=0.04$  \\
		\hline
	& 0.004		& 0.149			& 0.142		& 0.235	&	0.394	\\
4	& 0.006		& 0.085		& 0.156		& 0.257	&	0.314	\\
	& 0.008		& 0.086		& 0.148		& 0.231		&	0.326  \\
	\hline
	& 0.004		& 0.106			& 0.168		& 0.224	&	0.259	\\
5	& 0.006		& 0.082		& 0.149		& 0.216		&	0.272	\\
	& 0.008		& 0.084		& 0.130		& 0.196		&	0.287	\\
	\hline
	& 0.004		& 0.128			& 0.227	& 0.220	&	0.252	\\
6	& 0.006		& 0.099		& 0.176		& 0.213		&	0.289	\\
	& 0.008		& 0.086		& 0.141		& 0.216		&	0.265	\\
	\hline
	& 0.004		& 0.211			& 0.135		& 0.321		&	0.265	\\
7	& 0.006		& 0.086		& 0.146		& 0.259 	&	0.288	\\
	& 0.008		& 0.095		& 0.130		& 0.217		&	0.302	\\
	\hline
	& 0.004		& 0.101			& 0.151		& 0.233	&	0.246	\\
8	& 0.006		& 0.086		& 0.123		& 0.170		&	0.234	\\
	& 0.008		& 0.089		& 0.144		& 0.218		&	0.251	\\
	\hline
	& 0.004		& 0.101			& 0.145		& 0.194		&	0.249	\\
9	& 0.006		& 0.085		& 0.134		& 0.166		&	0.226	\\
	& 0.008		& 0.090		& 0.135		& 0.205	&	0.238	\\
	\hline
	& 0.004		& 0.131			& 0.139		& 0.205	&	0.222	\\
10	& 0.006		& 0.081			& 0.119		& 0.177		&	0.223	\\
	& 0.008		& 0.090		& 0.133		& 0.187		&	0.224	\\
	\hline
\end{tabular}
\end{center}
\caption{Dependence of $\varepsilon$ on the number of measurements $M$ and the parameter $\gamma$ for reconstruction of circles with noisy data $\epsilon = 0.01$, $0.02$, $0.03$ and $0.04$}
\label{tab:circles_err}
\end{table}

\section{Conclusion}

This paper presents the numerical algorithm for the reconstruction of the piecewise constant  coefficient. The performance of the algorithm is validated on numerical experiments. We have determined that the number of measurements is a crucial parameter for achieving good accuracy.  The smoothing parameter  $\gamma$ also affects the reconstruction results. The optimal values of these parameters can be determined experimentally. 

In this work,  we assumed that an object to be imaged contains only two materials with known piecewise constant coefficients. More general cases with several unknown coefficients will be a subject of our future work. It is well known that the complete electrode model is better suited for electrical impedance tomography than the continuum model considered here. We will extend the proposed algorithm to work with the complete electrode model in the future.


 

\begin{thebibliography}{10}

\bibitem{Chambers1999}
J.~Chambers, R.~Ogilvy, P.~I. Meldrum, and J.~Nissen.
\newblock {3D resistivity imaging of buried oil- and tar-contaminated waste
  deposits}.
\newblock {\em European Journal of Environmental and Engineering Geophysics},
  4(1):3--14, 1999.

\bibitem{Chambers2006}
J.~E. Chambers, O.~Kuras, P.~I. Meldrum, R.~D. Ogilvy, and J.~Hollands.
\newblock {Electrical resistivity tomography applied to geologic,
  hydrogeologic, and engineering investigations at a former waste-disposal
  site}.
\newblock {\em Geophysics}, 71(6), 2006.

\bibitem{Bodenstein2009a}
M.~Bodenstein, M.~David, and K.~Markstaller.
\newblock {Principles of electrical impedance tomography and its clinical
  application}.
\newblock {\em Critical Care Medicine}, 37(2):713--724, 2009.

\bibitem{Putensen2019}
C.~Putensen, B.~Hentze, S.~Muenster, and T.~Muders.
\newblock {Electrical Impedance Tomography for Cardio-Pulmonary Monitoring}.
\newblock {\em Journal of Clinical Medicine}, 8(8):1176, 2019.

\bibitem{McCarter1989}
W.~J. McCarter and S.~Garvin.
\newblock {Dependence of electrical impedance of cement-based materials on
  their moisture condition}.
\newblock {\em Journal of Physics D: Applied Physics}, 22(11):1773--1776, 1989.

\bibitem{Karhunen2010}
K.~Karhunen, A.~Sepp{\"{a}}nen, A.~Lehikoinen, P.~J.M. Monteiro, and J.~P.
  Kaipio.
\newblock {Electrical Resistance Tomography imaging of concrete}.
\newblock {\em Cement and Concrete Research}, 40(1):137--145, 2010.

\bibitem{Calderon1980}
A.~P. Calderon.
\newblock {On an inverse boundary value problem}.
\newblock In {\em Seminar on Numerical Analysis and its applications to
  Continuum Physics (Soc. Brasileira de Mat`ematica, Rio de Janeiro)}, pages
  65--73, 1980.

\bibitem{Sylvester1986}
J.~Sylvester and G.~Uhlmann.
\newblock {A uniqueness theorem for an inverse boundary value problem in
  electrical prospection}.
\newblock {\em Communications on Pure and Applied Mathematics}, 39(1):91--112,
  1986.

\bibitem{Novikov1988}
R.~G. Novikov.
\newblock {Multidimensional inverse spectral problem for the equation -$\Delta\psi$ + (v(x) - Eu(x))$\psi$ = 0}.
\newblock {\em Functional Analysis and Its Applications}, 22(4):263--272, 1988.

\bibitem{Nachman1996}
A.~Nachman.
\newblock {Global uniqueness for a two-dimensional inverse boundary value
  problem}.
\newblock {\em Annals of Mathematics}, 143(1):71--96, 1996.

\bibitem{Brown1997}
R.~M. Brown and G.~A. Uhlmann.
\newblock {Uniqueness in the inverse conductivity problem for nonsmooth
  conductivities in two dimensions}.
\newblock {\em Communications in Partial Differential Equations},
  22(5-6):1009--1027, 1997.

\bibitem{Astala2006}
K.~Astala and L.~P{\"{a}}iv{\"{a}}rinta.
\newblock {Calder{\'{o}}n's inverse conductivity problem in the plane}.
\newblock {\em Annals of Mathematics}, 163(1):265--299, 2006.

\bibitem{Alessandrini1988}
G.~Alessandrini.
\newblock {Stable determination of conductivity by boundary measurements}.
\newblock {\em Applicable Analysis}, 27:153--172, 1988.

\bibitem{Allers1991}
A.~Allers and F.~Santosa.
\newblock {Stability and resolution analysis of a linearized problem in
  electrical impedance tomography}.
\newblock {\em Inverse Problems}, 7(4):515--533, 1991.

\bibitem{Barcelo2001a}
J.~A. Barcel{\'{o}}, T.~Barcel{\'{o}}, and A.~Ruiz.
\newblock {Stability of the Inverse Conductivity Problem in the Plane for less
  regular conductivities}.
\newblock {\em Journal of Differential Equations}, 173:231--270, 2001.

\bibitem{Cheney1999}
M.~Cheney, D.~Isaacson, and J.~C. Newell.
\newblock {Electrical impedance tomography}.
\newblock {\em SIAM Review}, 41(1):85--101, 1999.

\bibitem{Borcea2002}
L.~Borcea.
\newblock {Electrical impedance tomography}.
\newblock {\em Inverse Problems}, 18:R99--R136, 2002.

\bibitem{Bruhl2000}
M.~Br{\"{u}}hl and M.~Hanke.
\newblock {Numerical implementation of two noniterative methods for locating
  inclusions by impedance tomography}.
\newblock {\em Inverse Problems}, 16(4):1029--1042, 2000.

\bibitem{Gebauer2007}
B.~Gebauer and N.~Hyv{\"{o}}nen.
\newblock {Factorization method and irregular inclusions in electrical
  impedance tomography}.
\newblock {\em Inverse Problems}, 23(5):2159--2170, 2007.

\bibitem{Somersalo1991}
E.~Somersalo, M.~Cheney, D.~Isaacson, and E.~Isaacson.
\newblock {Layer stripping: A direct numerical method for impedance imaging}.
\newblock {\em Inverse Problems}, 7(6):899--926, 1991.

\bibitem{Sylvester1992}
J.~Sylvester.
\newblock {A Convergent layer stripping algorithm for the radially symmetric
  impedence tomography problem}.
\newblock {\em Communications in Partial Differential Equations},
  17(11-12):1955--1994, 1992.

\bibitem{Siltanen2001}
S.~Siltanen, J.~Mueller, and D.~Isaacson.
\newblock {An implementation of the reconstruction algorithm of a Nachman for
  the 2D inverse conductivity problem}.
\newblock {\em Inverse Problems}, 16(3):681--699, 2000.

\bibitem{Knudsen2009}
K.~Knudsen, M.~Lassas, J.~L. Mueller, and S.~Siltanen.
\newblock {Regularized d-bar method for the inverse conductivity problem}.
\newblock {\em Inverse Problems and Imaging}, 3(4):599--624, 2009.

\bibitem{Cheney1990}
M.~Cheney, D.~Isaacson, J.~C. Newell, S.~Simske, and J.~Goble.
\newblock {NOSER: An algorithm for solving the inverse conductivity problem}.
\newblock {\em International Journal of Imaging Systems and Technology},
  2(2):66--75, 1990.

\bibitem{Bruhl2001}
M.~Br{\"{u}}hl.
\newblock {Explicit characterization of inclusions in electrical impedance
  tomography}.
\newblock {\em SIAM Journal of Mathematical Analysis}, 32(6):1327--1341, 2001.

\bibitem{Yorkey1987}
T.~J. Yorkey, J.~G. Webster, and W.~J. Tompkins.
\newblock {Comparing Reconstruction Algorithms for Electrical Impedance
  Tomography}.
\newblock {\em IEEE Transactions on Biomedical Engineering},
  BME-34(11):843--852, 1987.

\bibitem{Dobson1992}
D.~C. Dobson.
\newblock {Convergence of a reconstruction method for the inverse conductivity
  problem}.
\newblock {\em SIAM Journal on Applied Mathematics}, 52(2):442--458, 1992.

\bibitem{Borcea2001}
L.~Borcea.
\newblock {A nonlinear multigrid for imaging electrical conductivity and
  permittivity at low frequency}.
\newblock {\em Inverse Problems}, 17(2):329--359, 2001.

\bibitem{Wexler1985}
A.~Wexler, B.~Fry, and M.~R. Neuman.
\newblock {Impedance-computed tomography algorithm and system}.
\newblock {\em Applied Optics}, 24(23):3985, 1985.

\bibitem{Kohn1990}
R.~V. Kohn and A.~McKenney.
\newblock {Numerical implementation of a variational method for electrical
  impedance tomography}.
\newblock {\em Inverse Problems}, 6(3):389--414, 1990.

\bibitem{Dorn2000}
O.~Dorn, E.~. Miller, and C.~M. Rappaport.
\newblock A shape reconstruction method for electromagnetic tomography using
  adjoint fields and level sets.
\newblock {\em Inverse Problems}, 16(5):1119--1156, 2000.

\bibitem{Osher1988}
S.~Osher and J.~A. Sethian.
\newblock {Fronts propagating with curvature-dependent speed: Algorithms based
  on Hamilton-Jacobi formulations}.
\newblock {\em Journal of Computational Physics}, 79(1):12--49, 1988.

\bibitem{Santosa1996}
F.~Santosa.
\newblock {A Level-set Approach for Inverse Problems Involving Obstacles}.
\newblock {\em Esaim : Control, Optimisation and Calculus of Variations},
  1:17--33, 1996.

\bibitem{Osher2001}
S.~J. Osher and F.~Santosa.
\newblock {Level Set Methods for Optimization Problems Involving Geometry and
  Constraints I. Frequencies of a Two-Density Inhomogeneous Drum}.
\newblock {\em Journal of Computational Physics}, 171(1):272--288, 2001.

\bibitem{Ito2001}
K.~Ito, K.~Kunisch, and Z.~Li.
\newblock {Level-set function approach to an inverse interface problem}.
\newblock {\em Inverse Problems}, 17(5):1225--1242, 2001.

\bibitem{Chan2004}
T.~F. Chan and X.~C. Tai.
\newblock {Level set and total variation regularization for elliptic inverse
  problems with discontinuous coefficients}.
\newblock {\em Journal of Computational Physics}, 193(1):40--66, 2004.

\bibitem{Chung2005}
E.~T. Chung, T.~F. Chan, and X.~C. Tai.
\newblock {Electrical impedance tomography using level set representation and
  total variational regularization}.
\newblock {\em Journal of Computational Physics}, 205(1):357--372, 2005.

\bibitem{Agnelli2018}
J.~P. Agnelli, A.~{De Cezaro}, and A.~Leit{\~{a}}o.
\newblock {A regularization method based on level sets and augmented Lagrangian
  for parameter identification problems with piecewise constant solutions}.
\newblock {\em Inverse Problems}, 34(12):0--16, 2018.

\bibitem{Lin2018}
G.~Lin, Y.~Cheng, and Y.~Zhang.
\newblock {A parametric level set based collage method for an inverse problem
  in elliptic partial differential equations}.
\newblock {\em Journal of Computational and Applied Mathematics}, 340:101--121,
  2018.

\bibitem{Liu2019}
D.~Liu, D.~Smyl, and J.~Du.
\newblock {A parametric level set-based approach to difference imaging in
  electrical impedance tomography}.
\newblock {\em IEEE Transactions on Medical Imaging}, 38(1):145--155, 2019.

\bibitem{Tai2004}
X.~C. Tai and T.~F. Chan.
\newblock {A survey on multiple level set methods with applications for
  identifying piecewise constant functions}.
\newblock {\em International Journal of Numerical Analysis and Modeling},
  1(1):25--47, 2004.

\bibitem{Burger2005a}
M.~Burger and S.~J. Osher.
\newblock {A survey on level set methods for inverse problems and optimal
  design}.
\newblock {\em European Journal of Applied Mathematics}, 16(2):263--301, 2005.

\bibitem{Gibou2018}
F.~Gibou, R.~Fedkiw, and S.~Osher.
\newblock {A review of level-set methods and some recent applications}.
\newblock {\em Journal of Computational Physics}, 353:82--109, 2018.

\bibitem{Kolesov2019}
A.~E. Kolesov, D.~Kh. Ivanov, and P.~N. Vabishchevich.
\newblock {Recovery of a piecewise constant lower coefficient of an elliptic
  equation}.
\newblock {\em Journal of Physics: Conference Series}, 1392(1), 2019.

\bibitem{Ivanov2019}
D.~Kh. Ivanov, A.~E. Kolesov, and P.~N. Vabishchevich.
\newblock {Numerical method for recovering the piecewise constant right-hand
  side function of an elliptic equation from a boundary overdetermination
  data}.
\newblock {\em Journal of Physics: Conference Series}, 1392(1):0--6, 2019.

\bibitem{Dorn2006}
O.~Dorn and D.~Lesselier.
\newblock {Level set methods for inverse scattering}.
\newblock {\em Inverse Problems}, 22(4):R67--R131, 2006.

\bibitem{AlnaesBlechta2015a}
M.~S. Aln{\ae}s, J.~Blechta, J.~Hake, A.~Johansson, B.~Kehlet, A.~Logg,
  C.~Richardson, J.~Ring, M.~E. Rognes, and G.~N. Wells.
\newblock The FEniCS project version 1.5.
\newblock {\em Archive of Numerical Software}, 3(100), 2015.

\end{thebibliography}

\end{document}